\documentclass [12pt] {report}
\usepackage[cp1251]{inputenc}
\usepackage[dvips]{graphicx}
\usepackage{amssymb}
\newcommand {\Cos} {\mathop{\rm Cos}\nolimits}

\newcommand {\Sin} {\mathop{\rm Sin}\nolimits}
\newcommand {\e} {\mathop{\rm e}\nolimits}

\newcommand {\D}[2] {\displaystyle\frac{\partial{#1}}{\partial{#2}}}
\newcommand {\DD}[2] {\displaystyle\frac{\partial^2{#1}}{\partial{#2^2}}}
\newcommand {\Dd}[3] {\displaystyle\frac{\partial^2{#1}}{\partial{#2}\partial{#3}}}

\newcommand {\al} {\alpha}

\newcommand {\si} {\sigma}

\newcommand {\de} {\delta}
\newcommand {\De} {\Delta}

\newcommand {\iy} {\infty}
\newcommand {\prtl} {\partial}
\newcommand {\fr} {\displaystyle\frac}
\newcommand {\suml} {\sum\limits}

\newcommand {\be} {\begin{equation}}
\newcommand {\ee} {\end{equation}}
\newcommand {\ba} {\begin{array}}
\newcommand {\ea} {\end{array}}
\newcommand {\bp} {\begin{picture}}
\newcommand {\ep} {\end{picture}}
\newcommand {\bc} {\begin{center}}
\newcommand {\ec} {\end{center}}

\newcommand {\bt} {\begin{tabular}}
\newcommand {\et} {\end{tabular}}
\newcommand {\lf} {\left}
\newcommand {\rg} {\right}

\newcommand {\nin}{\noindent}
\newcommand {\cB} {{\cal B}}
\newcommand {\cC} {{\cal C}}

\newcommand {\cE} {{\cal E}}
\newcommand {\cF} {{\cal F}}

\newcommand {\cI} {{\cal I}}

\newcommand {\cP} {{\cal P}}
\newcommand {\cQ} {{\cal Q}}
\newcommand {\cR} {{\cal R}}
\newcommand {\cS} {{\cal S}}

\newcommand {\bP} {{\bf P}}

\newcommand {\bR} {{\bf R}}

\newcommand {\g}  {\stackrel{g\to -g}{\Longleftrightarrow}}

\newcommand {\ses} {\medskip}
\newcommand {\pgbrk}{\pagebreak}

\def\2#1#2#3{{#1}_{#2}\hspace{0pt}^{#3}}
\def\3#1#2#3#4{{#1}_{#2}\hspace{0pt}^{#3}\hspace{0pt}_{#4}}
\newcounter{sctn}
\def\sec#1.#2\par{\setcounter{sctn}{#1}\setcounter{equation}{0}
                  \noindent{\bf\boldmath#1.#2}\bigskip\par}
\begin {document}

\begin {titlepage}

\vspace{0.1in}

\begin{center}
{\Large
\bf
Finsleroid space with angle and scalar product
}\\
\end{center}

\vspace{0.3in}

\begin{center}

\vspace{.15in}
{\large G. S. Asanov\\}
\vspace{.25in}
{\it Division of Theoretical Physics, Moscow State University\\
119992 Moscow, Russia\\
(e-mail: asanov@newmail.ru)}
\vspace{.05in}

\end{center}

\begin{abstract}

A systematic approach has been developed to encompass the
Minkowski-type extension of Euclidean geometry such that a
one-vector  anisotropy is permitted, retaining simultaneously the
concept of angle. For the respective geometry, the Euclidean unit
ball is to be replaced by  the body which is  convex and rotund
and is found on assuming that its surface (the indicatrix
extending the unit sphere) is a space of constant positive
curvature. We have called the body the Finsleroid in view of  its
intrinsic relationship with the metric function of  Finsler
type.
 The main  point of the present paper is the angle coming
from geodesics through the cosine theorem,  the underlying idea
being to derive the angular measure from the solutions to the
geodesic equation  which prove to be obtainable in  simple
explicit  forms. The substantive items concern geodesics, angle,
scalar product, perpendicularity, and two-vector metric tensor.
The Finsleroid-space-associated one-vector Finslerian metric
function admits in quite a natural way an attractive two-vector
extension.

\end{abstract}

\end{titlepage}

\vskip 1cm

\bc
{\bf  0. Introduction}
\ec
\bigskip

The Euclidean geometry is  simple and totally spherically
symmetric, and corresponds well to our ordinary everyday
experience and intuition, while the Finsler or Banach-Minkowski
geometries [1-9] are much more extended and sophisticated constructions
that may serve to reflect various anisotropic scenarios. When a
single vector is distinguished geometrically   to be the only isotropic direction in
  extending the Euclidean geometry,
the sphere may not be regarded as an exact carrier of the
unit-vector image. So under respective conditions one may expect
that some directionally-anisotropic figure should be substituted
with the sphere. To this end we shall use the Finsleroid which,
being convex and rotund, is not, however, a second-order figure.
The constant positive curvature is the fundamental property of the
Finsleroid.

The present paper develops and elaborates in much detail the
related Finsleroid-geometry (initiated by the author earlier in
[10-12]) in the direction of evidencing the concepts of angle and scalar product.
 No special knowledge of Banach-Minkowski of Finsler geometries is assumed.

It will be recollected that, despite the fact that in geometry one
certainly  needs to use not only length but also angle and  scalar product,
various known   attempts  to introduce the concept of angle in the
Minkowski or Finsler spaces  were steadily encountered with
drawback positions:

\ses

"Therefore no particular angular measure can be entirely natural in
Minkowski geometry. This is evidenced by the innumerable attempts to define
such a measure, none of which found general acceptance``.
(Busemann [2], p. 279.)

\ses

"Unfortunately, there exists a number of distinct invariants
in a Minkowskian space all of which
reduce to the same classical
euclidean invariant if
the Minkowskian space
degenerates into
a euclidean space.
Consequently,
distinct definitions of
the trigonometric functions and of  angles
have appeared in the literature concerning
 Minkowskian and
Finsler spaces``. (Rund [3], p. 26)

\ses

A short but profound review of the respective attempts can be
found in Section 1.7 of the book [4]. The fact that the attempts
have never been unambiguous seems to be due to a lack of the
proper tools. For the opinion was taken for granted that the angle
ought to be defined or constructed in terms of the basic
Finslerian metric tensor (and whence ought to be explicated from
the initial Finslerian metric function). Let us doubt  the opinion
from the very beginning. Instead, we would like to raise
alternatively the principle that the angle
 is a concomitant of the geodesics (and  not  of the metric function proper).
The angle is determined by two vectors (instead of one vector in
case of the length) and actually implies using a due extension of
the Finslerian metric function to a two-vector metric function (to
a scalar product). Below, the principle is applying to the
Finsleroid space in a systematic way.

\ses

The abbreviations FMF and FMT will be used for the
Finsleroid metric function and the associated Finslerian metric
tensor, respectively. The notation  $\cE^{PD}_g$ will be applied to the Finsleroid space,
 with the upperscripts $``PD"$
meaning ``positive-definite".
The characteristic parameter $g$ may take on the values between $-2$ and $2$; at $g=0$
 the space is reduced to become an ordinary Euclidean one.

\ses
\ses

\clearpage

 \bc
{\bf Chapter 1: Synopsis of new conclusions}
\ec

\setcounter{sctn}{1}
\setcounter{equation}{0}
\ses
\ses

The angle $\al$ obtained in the Finsleroid Geometry under study
has the following remarkable property: if the consideration is
restricted to the $(N=2)$-dimensional Finsleroid-Minkowski plane,
then $\al$ {\it is a factor of the respective Landsberg angle}
(see  Section 1.1).
The respective  $\cE_g^{PD}$-Generalized  Trigonometric Functions
are appeared. Section 1.2 is devoted to reviewing the key
and basic concepts determined by the angle.
 Chapter 1 ends with Section 1.3 in which (on  returning the treatment from the
auxiliary quasi-Euclidean framework back to the primary Finsleroid
space) we are able to display the form of the associated two-vector metric tensor.

\ses \ses

1.1. {\it Finsleroid-Minkowski plane}. When reducing the
consideration to the {\it Minkowski plane}, with the dimension
$N=2$ and the orthogonalized form $r_{pq}=\de_{pq}$ of the input
Euclidean metric tensor, the Finsleroid-adapted vector components
(cf. the representation (2.83) in Chapter 2)   take on the form
\be
R^1=\fr K{hJ}\sin f,
\quad
 R^2=\fr KJ(\cos f-\fr12G\sin f),
 \ee
from which it follows that \be
R^2\D{R^1}f-R^1\D{R^2}f=\fr1{hJ^2}K^2. \ee
 Since also
$\sqrt{\det(g_{pq})}=J^2$ (cf. Eq. (2.64) in Chapter 2), from the
equality (1.2) we conclude
$$d\al_{Finsleroid-Landsberg}=\fr1h df
$$
 (see, e.g., p. 85 of [8] for the definition of the Landsberg
angle),
 where
$h$ is the constant (2.13) of Chapter 2.

Therefore, the following theorem is valid.

\ses \ses

 {\large Theorem} 1.1.
{\it Restricting the Finsleroid geometry to the Minkowski plane,
the quantity $f$ in the representation \rm  (1.1) \it is the
factor $h$ of the Landsberg angle. }

\ses \ses

It is also possible to draw

\ses \ses

{\large  Theorem} 1.2. {\it The Finsleroid Indicatrix on the
Minkowski plane is strongly convex.}

\ses \ses

{\it Proof.} Let us verify the relevant criterion formulated on p.
88 of [8]. In terms of our notation, we calculate accordingly:
$$
\fr {\DD{R^2}{f} \D{R^1}{f} - \D{R^2}{f} \DD{R^1}{f} } {\D{R^2}{f}
R^1 - R^2\D{R^1}{f} } = \fr{-\fr1{h^3}}{-\fr1h}=\fr1{h^2}.
$$
Since the right-hand side here is always positive, the criterion
works fine and, therefore,  Theorem 1.2 is valid.

\ses \ses

Noting that $ ds :=\sqrt{g_{pq}(g;R)dR^pdR^q}=\fr1hdf, $
we conclude that
\be
 ds=\fr1hdf.
  \ee
   In particular, the latter
equality entails

\ses\ses

{\large  Theorem} 1.3. {\it The  length $L_{\cI} :=\int ds$ of the
Finsleroid Indicatrix is \be L_{\cI}=\fr{2\pi}h\ge 2\pi,
 \ee
showing the properties }
\be L_{\cI}=2\pi \quad \mbox{\rm if and
only if}\quad g=0\quad \mbox{(the Euclidean case)}
\ee {\it and}
 \be L_{\cI}\to
\infty \quad {\rm when} \quad |g|\to 2. \ee

\ses \ses

From (1.1) and (1.3) it can readily be explicated that the {\it
Rund equation}
 \be
\fr{d^2R^p}{ds^2}+I\fr{dR^p}{ds}+R^p=0
\ee
 holds fine with
 \be
I=-g.
\ee If the meaning of the Cartan scalar is acquired to the
quantity $I$ thus appeared in (1.7) (cf. [8]), one may state the
following:

\ses \ses

{\large  Theorem} 1.4. {\it The Cartan scalar for the
Finsleroid-Minkowski plane is the constant which equals the
negative of the characteristic Finsleroid parameter $g$. }

\ses\ses

 Eqs. (1.1) suggest  naturally to propose the following
 {\it $\cE_g^{PD}$-Generalized  Trigonometric Functions}:
\ses\\
\be
\Cos f~:=\fr1J(\cos f-\fr G2\sin f), \qquad \Sin f~:=\fr1{hJ}\sin f,
\ee
and
\be
\Cos^*f~:=\fr1{h^2J}(\cos f+\fr G2\sin f).
\ee
They reveal the properties
\be
R^1=K\Sin f, \qquad  R^2=K\Cos f,
\ee
and
\be
(\Cos f)'=-\fr1h\Sin f, \qquad
(\Sin f)'=\Cos^*f,
\ee
together with
\be
(\Cos^* f)'=-\Sin f +G\Cos f,
\ee
where the prime stands for the derivative with respect to $f$.

\ses
\ses

\setcounter{sctn}{2}
\setcounter{equation}{0}

1.2. {\it
 Finsleroid angle.} Given two vectors $R_1\in V_N$ and $R_2\in
 V_N$.
Applying the quasi-Euclidean transformation (see (5.11) in Chapter
2) to Eq. (1.36) of Chapter 3 at any dimension $N\ge 2$ results in
the following {\it $\cE_g^{PD}$-scalar product}:
 \be
<R_1,R_2>=K(g;R_1)K(g;R_2) \cos \Bigl[ \fr1h\arccos\fr {
A(g;R_1)A(g;R_2)+h^2r_{be}R_1^bR_2^e } {
\sqrt{B(g;R_1)}\,\sqrt{B(g;R_2)} } \Bigl],
 \ee
 \ses \ses \ses so
that the {\it $\cE_g^{PD}$-angle}
 \be
 \al(R_1,R_2) = \fr1h\arccos \fr
{ A(g;R_1)A(g;R_2)+h^2r_{be}R_1^bR_2^e } {
\sqrt{B(g;R_1)}\,\sqrt{B(g;R_2)} }
 \ee
 is appeared between the vectors $R_1$
and $R_2$; the functions $B, K$, as well as the function  $A$ can
be found in Section 2.2 of Chapter 2.

In the Euclidean limit proper, the angle (2.2) is reduced to read merely
$$
 \al(R_1,R_2)_{{\Bigl |\Bigr.}_{g=0}} = \arccos \fr
{ R_1^NR_2^N+r_{be}R_1^bR_2^e } {
\sqrt{(R^N_1)^2+r_{be}R_1^bR_1^e}\,\sqrt{(R_2^N)^2+r_{be}R_2^bR_2^e} }.
$$

For the intermediate angle $\nu$ defined by Eq. (1.44) in Chapter
3, we obtain \be \nu=\arctan\fr{sK(g;R_2)\sin\al}{K(g;R_1)\De
s+[K(g;R_2)\cos\al-K(g;R_1)]s}. \ee

The
$\cE^{PD}_g$-space general solution to the geodesic equation (presented by
 Eqs. (2.88)-(2.89) in  Chapter 2) reads
\be R^p(s)=\mu^p(g;t(g;s)), \ee where $ t(g;s)$ is given by Eq.
(1.41) of Chapter 3, and $ \mu^p $ are the functions which realize
the quasi-Euclidean transformation according to Eqs. (5.14)-(5.15)
of  Chapter 2; $s$ is the arc-length parameter (see Eq. (2.89) in
Chapter 2). The relevant explicit formulas are


 \be R^N(s)=(t^N(s)-\fr12Gm(s))/k(s),
\qquad R^a(s)=\fr1ht^a(s)/k(s)\ee
 with
 \be
t^N(s)=\fr{K_s}{\sin(h\al)}\Bigl[\fr{A_1}{\sqrt{B_1}}\sin(h(\al-\nu))
+ \fr{A_2}{\sqrt{B_2}}\sin(h\nu)\Bigr], \ee
where
 \be
t^a(s)=h\fr{K_s}{\sin(h\al)}\Bigl[\fr{R_1^a}{\sqrt{B_1}}\sin(h(\al-\nu))
+ \fr{R_2^a}{\sqrt{B_2}}\sin(h\nu)\Bigr] ,\ee
 and
  \be
K_s=\sqrt{K^2(g;R_2)+2sK(g;R_1)\sqrt{1-\lf(\fr{K(g;R_2)\sin\al}{\De
s }\rg)^2}+s^2},\ee \ses
\be
k(s)=\exp\lf(\fr12G\arctan\fr{t^N(s)}{m(s)}\rg),
\ee
\ses
\be
m(s)=\sqrt{r_{ab}t^a(s)t^b(s)}.\ee

Along the geodesics, the behaviour law for the squared FMF is
quadratic:
\be
K^2(g;R(s))=a^2+2bs+s^2
 \ee
(cf. Eq. (1.12) in Chapter 3); $a$ and $b$ are integration constants.

\ses

Below the picture symbolizes the role which the angles (2.2) and (2.3)
are playing  in the geodesic line $C$ which joins
two points $P_1$ and $P_2$.

\ses
\ses

\begin{figure}[ht]
 \centering
\includegraphics[width=8cm]{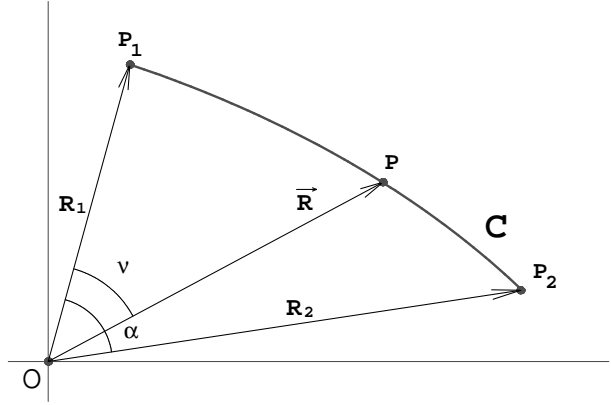}
\caption{[The geodesic $C$ and the angles $\al=\angle P_1OP_2$ and $\nu=\angle P_1OP$]}
\end{figure}

\ses
\ses

On this way  the following substantive items can be arrived at.

\ses
\ses

\nin
{\it
The
$\cE_g^{PD}$-Case
Cosine Theorem
}
\ses\\
\be (\Delta s)^2 = (K(g;R_1))^2 + (K(g;R_2))^2 - 2
K(g;R_1)K(g;R_2) \cos \al.
 \ee

\ses
\ses

\nin
{\it
The
$\cE_g^{PD}$-Case
Two-Point Length
}
\ses\\
\be | R_1 \ominus R_2 |^2 = (K(g;R_1))^2 + (K(g;R_2))^2 - 2
K(g;R_1)K(g;R_2) \cos \al. \ee

\ses
\ses

\nin
{\it
The
$\cE_g^{PD}$-Case
Scalar Product
}
\ses\\
\be <R_1,R_2>=K(g;R_1)K(g;R_2) \cos \al. \ee

\ses

\nin
{\it
At equal vectors, the reduction}
\be
<R,R>=
K^2(g;R)
\ee
takes place, that is, the two-vector scalar product
(2.1) reduces  exactly  to the squared FMF.

\ses \ses

\nin {\it The $\cE_g^{PD}$-Case Perpendicularity }
\ses\\
\be
<R,R^{\perp}>=K(g;R)K(g;R^{\perp}),
\ee
in which case $\al=\pi/2$.

\ses
\ses

\nin
{\it
The
$\cE_g^{PD}$-Case
Pythagoras Theorem
}
\ses\\
\be | R_1 \ominus R_2 |^2 = (K(g;R_1))^2 + (K(g;R_2))^2 \ee holds
fine.

\ses
\ses

\nin
{\it
The Finsleroid angle $\al$ between two vectors ranges over}
$$
0\le\al\le
\al_{max},
$$
where
$$
\al_{max}=\fr1h\pi\ge\pi \quad {\rm with~equality~if~and~only~if}~g=0,
$$
so that
$$
\al_{max}
\mathop{\Longrightarrow}_{g\to 2}\infty.
$$

\ses
\ses
\ses
\ses

The identification \be | R_2 \ominus R_1 |^2 =(\Delta s)^2 \ee
yields another lucid representation \be | R_1 \ominus R_2 |^2 =
(K(g;R_1))^2 + (K(g;R_2))^2 - 2 K(g;R_1)K(g;R_2) \cos \al\,. \ee

We observe the symmetry
\be
| R_1 \ominus R_2 | = | R_2 \ominus
R_1 |.
\ee

Particularly, from (2.2) it directly ensues that  the value of the
angle  $\al$ formed by  a vector $R$ with the Finsleroid
$R^N$-axis is given by
 \be \al = \fr1h\arccos \fr { A(g;R) } {
\sqrt{B(g;R)}\ }, \ee and with $(N-1)$-dimensional equatorial
$\{{\bf R}\}$-plane of Finsleroid is prescribed as \be \al =
\fr1h\arccos \fr { L(g;R) } { \sqrt{B(g;R)}\ }; \ee here, $ L$ is
the function (2.36) of Chapter 2.

\ses
\ses

1.3. {\it Two-vector metric tensor.} Let $R\in V_N$ and $S\in V_N$
be two vectors. From (2.1) we obtain \ses
\ses\\
$$
\D{<R,S>}{R^p}=R_p
\fr{
<R,S>
}
{K^2(g;R)}
+hK(g;S)s_p(g;R,S)\sin
\al
$$
and
\ses
\ses\\
$$
\D{<R,S>}{S^q}=S_q
\fr{
<R,S>
}
{K^2(g;S)}
+hK(g;R)s_p(g;S,R)\sin
\al.
$$
\ses
\ses
For the associated
{\it
$\cE_g^{PD}$-two-vector metric tensor}
$$
G_{pq}(g;R,S)
:=
\Dd{<R,S>}{S^q}{R^p}
$$
we can find explicitly the representation
\ses\
$$
G_{pq}(g;R,S)
=
\lf(
\fr{R_p}
{
K(g;R)
}
\fr
{S_q}
{K(g;S)
}
-h^2s_p(g;R,S)s_q(g;S,R)
\rg)
\cos\al
$$
\ses
\ses
\ses
\ses
$$
+
h\Bigr[
\lf(
\fr{R_p}
{
K(g;R)
}
s_q(g;S,R)
+
\fr{S_q}
{
K(g;S)
}
s_p(g;R,S)
\rg)
+s_{pq}(g;R,S)
\Bigr]
\sin\al,
$$
\ses
where
\ses\\
$$
s_p(g;R,S)=
\fr{M_p(g;R,S)}
{
W(g;R,S)
}
\fr{K(g;R)}{B(g;R)}
$$
and
\ses\\
$$
s_{pq}(g;R,S)=
K(g;S)\D{s_p(g;R,S)}{S^q}
$$
with
\ses\\
$$
W(g;R,S)
=
\sqrt{
B(g;R)B(g;S)-
\Bigl[
A(g;R)A(g;S)+h^2r_{be}R^bS^e
\Bigl]^2
}
$$
and
\ses\\
$$
M_p(g;R,S)=B(g;R)
\sqrt{B(g;R)}\,\sqrt{B(g;S)}
\fr1{h^2}\fr{\partial}{\partial R^p}
\fr
{
A(g;R)A(g;S)+h^2r_{be}R^bS^e
}
{
\sqrt{B(g;R)}\,\sqrt{B(g;S)}
}.
$$
The latter vector  has the components
\ses
$$
h^2M_N(g;R,S)=
B(g;R)A(g;S)-
\Bigl[
A(g;R)A(g;S)+h^2r_{be}R^bS^e
\Bigl]
A(g;R)
$$
and
$$
h^2M_a(g;R,S)=B(g;R)
\lf(\fr12g
\fr{R^b}{q(R)}
A(g;S)+h^2S^b\rg)
r_{ab}
$$
\ses
\ses
$$
-
\Bigl[
A(g;R)A(g;S)+h^2r_{be}R^bS^e
\Bigl]
\lf(\fr12gR^N
+q(R)
\rg)
\fr{R^b}{q(R)}
r_{ab},
$$
\ses
\ses
\ses
\ses
which can be simplified to get
$$
M_N(g;R,S)=
q^2(R)A(g;S)-
r_{be}R^bS^e
A(g;R)
$$
and
$$
M_a(g;R,S)=
\lf(
-R^NR^b
A(g;S)+S^bB(g;R)
-
r_{ec}R^eS^c
\lf(q(R)+\fr12gR^N\rg)
\fr{R^b}{q(R)}
\rg)
r_{ab}.
$$

The identity
$$
M_p(g;R,S)R^p=0
$$
holds.

The symmetry
\ses
\ses
$$
G_{pq}(g;
R,S)=
G_{qp}(g;
S,R)
$$
\ses
and the Finslerian limit
$$
\lim_{
S^r\to
R^r
}
\Bigl\{
G_{pq}(g;
R,S)
\Bigl\}=
g_{pq}(g;R)
$$
(cf. Eq. (2.6) in Chapter 3) can straightforwardly be verified;
the components $ g_{pq}(g;R) $ are presented in Chapter 2 by the
list (2.60)-(2.61). The two-vector $\cE_g^{PD}$-metric tensor can
also be obtained as  the transform
$$
G_{pq}(g;
R,S)
=
\si^r_p(g;R)
\si^s_q(g;S)
n_{rs}(g;
{\bf t}_1,
{\bf t}_2)
$$
(cf. Eq. (5.47) in Chapter 2)
of the two-vector quasi-Euclidean tensor (see Eq. (2.2) in  Chapter 3),
where
$$
t^r_1=\si^r(g;R),\qquad
t^s_2=\si^s(g;S)
$$
(cf. Eqs. (5.10) and (5.35) in  Chapter 2).

\ses
\ses

\bc
{\bf Chapter 2. Centerpiece properties}
\ec

\setcounter{sctn}{1}
\setcounter{equation}{0}

2.1. {\it  Motivation}. Below Section 2.2 gives an account of the
notation and conventions for the space
 $\cE^{PD}_g$
and introduces the initial concepts and definitions that are
required. The space is constructed by assuming an axial symmetry
and, therefore, incorporates a single preferred direction, which
we shall often refer as the $Z$-axis. After preliminary
introducing a characteristic quadratic form $B$, which is distinct
from the Euclidean sum of squares by entrance of a mixed term (see
Eq. (2.22)), we define the FMF $K$ for the space $\cE^{PD}_g$ by
the help of the formulae (2.30)-(2.33). A characteristic feature
of the formulae is the occurrence of the function $``\arctan"$.
Next, we calculate basic tensor quantities of the space. There
appears a remarkable phenomenon, which essentially simplifies all
the constructions, that the associated Cartan tensor occurs being
of a simple algebraic structure (see Eqs. (2.66)-(2.67)). In
particular, the phenomenon gives rise to a simple structure of the
 associated curvature tensor (Eq. (2.69)).
As well as in the Euclidean geometry  the locus of the unit
vectors issuing from fixed  point of origin is the unit sphere, in
the $\cE_g^{PD}$-geometry under development the locus is the
boundary (surface) of the Finsleroid. We call the boundary the
Finsleroid Indicatrix. It can rigorously be proved that the
Finsleroid Indicatrix  is a closed, regular, and strongly convex
(hyper)surface.The value of the curvature depends on the parameter
$g$ according to the simple law (2.73). The determinant of the
associated FMT is strongly positive in accordance with Eqs.
(2.64)-(2.65).

In Section 2.3 we show how the consideration can conveniently be
converted  into the co-approach. We derive the associated
Finsleroid Hamiltonian metric function. The explicit form thereof
is entirely similar to the form of the FMF $K$ up to the
substitution of $-g$ with $g$. The symmetry between the relevant
input  representations (2.30)-(2.38) and (3.20)-(3.28) is
astonishingly perfect.

Section 2.4 gathers together the lucid facts concerning details of
the form of the Finsleroids and co-Finsleroids. The Finsleroid is
a generalization of the unit ball and may be visualized as
comprising a stretched surface of revolution.

Its form essentially depends on the value of the
characteristic parameter, $g$. Under changing the sign of $g$, the
Finsleroid turns up with respect to its equatorial section. When
$|g|\to 2$, the Finsleroid is extending ultimately tending in its
form more and more to the cone.
 The form and all the properties of
the co-Finsleroid are essentially similar to that of the
Finsleroid of the opposite sign of the parameter $g$. Various
Maple9-designed figures have been presented to elucidate patterns
and details, and to make this Finsleroid-framework a plausible one
in methodological as well  educational respects,
 --  which also show all the basic features of the Finsleroids.

The $\cE^{PD}_g$-space has an auxiliary quasi-Euclidean structure,
 which is deeply inherent in the development.
Section 2.5 introduces for  the $\cE^{PD}_g$-space the
quasi-Euclidean map under which the Finsleroid goes into the unit
ball. The quasi-Euclidean space is simple in many aspects,
 so that relevant transformations make reduce various calculations and may
provide one with constructive ideas.

Motivated by these observations, in Section 2.6 we offer the
nearest interesting properties of the quasi-Euclidean metric
tensor (which is not of a Finslerian type). The tensor is a linear
combination of the Euclidean metric tensor and the product of two
unit vectors. Therefore, the basic relevant geometric objects,
including the Christoffel symbols, the curvature tensor, the
orthonormal frames, and the Ricci rotation coefficients, are
calculated in simple forms.
 The quasi-Euclidean
space is not flat, but proves to be conformally-flat.

\ses
\ses


\setcounter{sctn}{2}
\setcounter{equation}{0}

2.2. {\it Bases}. Suppose we are given an $N$-dimensional centered
vector space $V_N$ with some point $``O"$ being the origin. Denote
by $R$ the vectors constituting the space, so that $R\in V_N$ and
it is assumed that $R$ is issued from the point $``O"$.
 Any given
vector $R$ assigns a particular direction in $V_N$. Let us fix a
member $R_{(N)}\in V_N$, introduce the straight line $e_N$ oriented
along the vector $R_{(N)}$,
 and use this
$e_N$
to serve as a $R^N$-coordinate axis
in $V_N$.
In this way we get the topological product
\be
V_N=V_{N-1}\times e_N
\ee
together with the separation
\be
R=\{\bR,R^N\}, \qquad R^N\in e_N \quad {\rm and} \quad \bR\in V_{N-1}.
\ee
For convenience, we shall frequently use the notation
\be
R^N=Z
\ee
and
\be
R=\{\bR,Z\}.
\ee
Also, we introduce a Euclidean metric
\be
q=q(\bR)
\ee
over the $(N-1)$-dimensional vector space
$V_{N-1}$.

With respect to an admissible coordinate basis $\{e_a\}$ in
$V_{N-1}$, we obtain the coordinate representations \be
\bR=\{R^a\}=\{R^1,\dots,R^{N-1}\} \ee and \be
R=\{R^p\}=\{R^a,R^N\}\equiv\{R^a,Z\}, \ee together with \be
q(\bR)=\sqrt{r_{ab}R^aR^b}; \ee the matrix $(r_{ab})$ is assumed
to be symmetric and positive-definite. The indices $(a,b,\dots)$
and $(p,q,\dots)$ will be specified over the ranges
$(1,\dots,N-1)$ and $(1,\dots,N)$, respectively; vector indices
are up, co-vector indices are down; repeated up-down indices are
automatically summed;  $\de$ will stand for the Kronecker symbol,
such that  $({\de}_{ab})={\rm diag}(1,1,\dots)$. The variables \be
w^a=R^a/Z, \quad w_a=r_{ab}w^b, \quad w= q/Z, \ee where \be
w\in(-\iy,\iy), \ee are convenient whenever $Z\ne0$. Sometimes we
shall mention the associated metric tensor \be
r_{pq}=\{r_{NN}=1,~r_{Na}=0,~r_{ab}\} \ee meaningful over the
whole vector space $V_N$.

Given a parameter $g$ subject to ranging
\be
-2<g<2,
\ee
we introduce the convenient notation
\be
h=\sqrt{1-\fr14g^2},
\ee
\ses
\be
G=g/h,
\ee
\ses
\be
g_+=\fr12g+h, \qquad g_-=\fr12g-h,
\ee
\medskip
\be
g^+=-\fr12g+h, \qquad g^-=-\fr12g-h,
\ee
so that
\be
 g_++g_-=g, \qquad g_+-g_-=2h,
\ee
\medskip
\be
 g^++g^-=-g, \qquad g^+-g^-=2h,
\ee
\ses
\be
(g_+)^2+(g_-)^2=2,
\ee
\ses
\be
(g^+)^2+(g^-)^2=2.
\ee
The symmetry
\be
g_+\g -g_-, \qquad g^+\g -g^-
\ee
holds.


The {\it  characteristic
quadratic form}
\be
B(g;R) :=Z^2+gqZ+q^2
\equiv\fr12\Bigl[(Z+g_+q)^2+(Z+g_-q)^2\Bigr]>0
\ee
is of the negative discriminant, namely
\be
D_{\{B\}}=-4h^2<0,
\ee
because of Eqs. (2.12) and (2.13).
Whenever $Z\ne0$, it is also convenient to use the quadratic form
\be
Q(g;w) := B/(Z)^2,
\ee
obtaining
\be
Q(g;w)=1+gw+w^2>0,
\ee
together with the function
\be
E(g;w) :=
1+\fr12gw.
\ee
The identity
\be
E^2+h^2w^2=Q
\ee
can readily be verified.
In the limit $g\to 0$,
the definition (2.22) degenerates to the
 quadratic form of the input metric tensor (2.11):
\be
B|_{_{g=0}}=r_{pq}R^pR^q.
\ee
Also
\be
Q|_{_{g=0}}=1+w^2.
\ee


In terms of this notation, we propose the FMF
\be
K(g;R)=
\sqrt{B(g;R)}\,J(g;R),
\ee
where
\be
J(g;R)=\e^{\frac12G\Phi(g;R)},
\ee
\medskip
\medskip
\be
\Phi(g;R)=
\fr{\pi}2+\arctan \fr G2-\arctan\Bigl(\fr{q}{hZ}+\fr G2\Bigr),
\qquad {\rm if} \quad Z\ge 0,
\ee
\medskip
\be
\Phi(g;R)=
-\fr{\pi}2+\arctan \fr G2-\arctan\Bigl(\fr{q}{hZ}+\fr G2\Bigr),
\qquad  {\rm if}  \quad Z\le 0,
\ee
\ses\\
or in other convenient forms,
\ses\\
\be
\Phi(g;R)=
\fr{\pi}2+\arctan \fr G2-\arctan\Bigl(\fr{L(g;R)}{hZ}\Bigr),
\qquad  {\rm if}  \quad Z\ge 0,
\ee
\medskip
\be \Phi(g;R)= -\fr{\pi}2+\arctan \fr
G2-\arctan\Bigl(\fr{L(g;R)}{hZ}\Bigr), \qquad  {\rm if}  \quad
Z\le 0, \ee where \be L(g;R) =q+\fr g2Z, \ee and also
\be
\Phi(g;R)= \arctan\fr{A(g;R)}{hq}, \ee
\ses\\
where
\be
A(g;R)=Z+\fr12gq.
\ee
\ses\\
This FMF has been normalized to show the  properties
\ses\\
\be
-\fr{\pi}2
\le\Phi\le
\fr{\pi}2
\ee
\medskip
and
\be
\Phi=
\fr{\pi}2,
\quad {\rm if} \quad q=0 \quad {\rm and} \quad Z>0;
\qquad
\Phi=
-\fr{\pi}2,\quad {\rm if} \quad q=0 \quad {\rm and} \quad Z<0.
\ee

We also have
\ses\\
\be
\tan\Phi=\fr A{hq}
\ee
\ses
and
\be
\Phi|_{_{Z=0}}=
\arctan \fr G2.
\ee
\ses\\
It is often convenient to use the sign indicator
$
\epsilon_Z
$
 for the argument $Z$:
\ses
\ses
\be
\epsilon_Z=1, \quad {\rm if}\quad Z>0;
\quad
\epsilon_Z=-1, \quad {\rm if}\quad Z<0;
\quad
\epsilon_Z=0, \quad {\rm if}\quad Z=0.
\ee

\ses

Under these conditions, we call the considered space
the
{\it
$\cE^{PD}_g$-space}:
\be
\cE^{PD}_g :=\{V_N=V_{N-1}\times e_N;\,R\in V_N;\,K(g;R);\,g\}.
\ee


The right-hand part of the definition (2.30)
can be considered to be a function
$\breve K$
of  the arguments
$\{g;q,Z\}$,
such that
\be
\breve K(g;q,Z)
=K(g;R).
\ee
We observe that
\be
\breve K(g;q,-Z)\ne \breve K(g;q,Z),\qquad unless \quad g=0.
\ee
Instead, the function $\breve K$ shows the property of
$gZ$-\it parity\rm
\be
\breve K(-g;q,-Z)=\breve K(g;q,Z).
\ee
The $(N-1)$-space reflection invariance
holds true
\be
 K(g;R)\stackrel{R^a\leftrightarrow -R^a}{\Longleftrightarrow} K(g;R)
\ee (such an operation does not influence  the quantity $q$).


It is frequently convenient to rewrite the representation (2.30) in the form
\be
K(g;R)=|Z|V(g;w),
\ee
whenever $Z\ne0$, with
the
{\it
generating metric
function}
\be
V(g;w)=\sqrt{Q(g;w)}\,
j(g;w).
\ee
We have
$$
j(g;w)=J(g;1,w).
$$
Using (2.25) and (2.31)-(2.35), we obtain
\be
V'=wV/Q,
\qquad
V''=V/Q^2,
\ee
\ses
\be
(V^2/Q)'=-gV^2/Q^2,  \qquad (V^2/Q^2)'=-2(g+w)V^2/Q^3,
\ee
\ses
\be
j'=-\fr12gj/Q,
\ee
\bigskip\\
and also
\be
\fr12(V^2)'=wV^2/Q,
\qquad\quad
\fr12(V^2)''=(Q-gw)V^2/Q^2,
\ee
\ses
\be
\fr14(V^2)'''=-gV^2/Q^3,
\ee
together with
\be
\Phi'=-h/Q,
\ee
where the prime ($'$) denotes the differentiation with respect to~$w$.

\ses

Also,
\be
(A(g;R))^2+h^2q^2=B(g;R)
\ee
and
\be
(L(g;R))^2+h^2Z^2=B(g;R).
\ee
\ses


The simple results for these derivatives reduce
the task of computing the components of the associated FMT
to an easy exercise, indeed:
$$
R_p :=\fr12\D{K^2(g;R)}{R^p}:
$$
\ses\ses
\be
R_a=
r_{ab}R^b
\fr{K^2}{B},
\qquad
R_N=(Z+gq)
\fr{K^2}{B};
\ee
\ses
\ses
\ses
\ses
$$
g_{pq}(g;R)
:=
\fr12\,
\fr{\prtl^2K^2(g;R)}{\prtl R^p\prtl R^q}
=\fr{\prtl R_p(g;R)}{\prtl R^q}:
$$
\ses
\ses
\ses
\ses
\be
 g_{NN}(g;R)=[(Z+gq)^2+q^2]
\fr{K^2}{B^2},
\qquad
g_{Na}(g;R)=gq
r_{ab}R^b
\fr{K^2}{B^2},
\ee
\ses
\ses
\be
g_{ab}(g;R)=
\fr{K^2}{B}
r_{ab}-g\fr{
r_{ad}R^d
r_{be}R^e
Z
}{q}
\fr{K^2}{B^2}.
\ee
\ses\\
The reciprocal tensor components are
\ses
\be
g^{NN}(g;R)=(Z^2+q^2)
\fr1{K^2},
\qquad
g^{Na}(g;R)=-gqR^a
\fr1{K^2},
\ee
\ses
\be
g^{ab}(g;R)=
\fr{B}{K^2}
r^{ab}+g(Z+gq)\fr{R^aR^b}{q}
\fr1{K^2}.
\ee
\ses\\
The determinant of the FMT
given by Eqs. (2.60)-(2.61) can readily be found in the form
\ses
\be
\det(g_{pq}(g;R))=[J(g;R)]^{2N}\det(r_{ab})
\ee
\ses\\
which shows, on noting (2.31)-(2.33), that
\ses
\be
\det(g_{pq})>0 {\it \quad over~all~the~space} \quad V_N.
\ee

The associated {\it angular metric tensor}
$$
h_{pq} := g_{pq}-R_pR_q\fr1{K^2}
$$
proves to be given by the components
$$
h_{NN}(g;R)=q^2
\fr{K^2}{B^2},
\qquad h_{Na}(g;R)=-Z
r_{ab}R^b
\fr{K^2}{B^2},
$$
\ses
$$
h_{ab}(g;R)=
\fr{K^2}{B}
r_{ab}-(gZ+q)\fr{
r_{ad}R^d
r_{be}R^e
}q
\fr{K^2}{B^2},
$$
\ses\\
which entails
$$
\det(h_{ab})=\det(g_{pq})\fr1{V^2}.
$$


The  use of the components of the Cartan tensor
(given explicitly in the end of the present section)
 leads,
after rather tedious straightforward calculations, to
the following simple and remarkable result.

\ses
\ses

{\large   Theorem} 2.1.
\it The Cartan tensor associated with the FMF
\rm(2.30) \it is of the following special
algebraic form:
\be
C_{pqr}=\fr1N\lf(h_{pq}C_r+h_{pr}C_q+h_{qr}C_p-\fr1{C_sC^s}C_pC_qC_r\rg)
\ee
with
\be
C_tC^t=\fr{N^2}{4K^2}g^2.
\ee
\rm
\ses

Elucidating the structure of
the respective curvature tensor
\be
S_{pqrs} := (C_{tqr}\3Cpts-C_{tqs}\3Cptr)
\ee
results in the simple representation
\be
S_{pqrs}=-\fr{C_tC^t}{N^2}(h_{pr}h_{qs}-h_{ps}h_{qr}).
\ee
Inserting here (2.67), we are led to
\ses

{\large   Theorem} 2.2.
\it The curvature tensor of the space
$\cE^{PD}_g$ is of the special type
\be
S_{pqrs}=S^*(h_{pr}h_{qs}-h_{ps}h_{qr})/K^2
\ee
with \rm
\be
S^*=-\fr14g^2.
\ee

\ses
\ses

{\large Definition}.\, FMF (2.30) generates the {\it Finsleroid}
 \be \cF^{PD}_g:=\{R\in V_N:K(g;R)\le 1\}. \ee

\ses \ses

{\large Definition}.\, The {\it Finsleroid Indicatrix} $
\cI^{PD}_g $ is the boundary of the Finsleroid:
 \be \cI^{PD}_g :=\{R\in V_N:K(g;R)=1\}.
  \ee

\ses \ses

{\large  Note}. Since at $g=0$ the space $ \cE^{PD}_g $ is
Euclidean, then the body $ \cF^{PD}_{g=0} $ is a unit ball and $
\cI^{PD}_{g=0} $ is a unit sphere.

\ses
\ses


Recalling the known formula $ \cR=1+S^* $ for the indicatrix
curvature (see Section 1.2 in [4]), from (2.71) we conclude that
\be \cR_{Finsleroid ~ Indicatrix }=h^2=1-\fr14g^2,\ee so that $$ 0
< \cR_{Finsleroid ~Indicatrix } \le 1 $$ and
$$
\cR_{Finsleroid~ Indicatrix }\stackrel{g\to 0}{\Longrightarrow}
\cR_{Euclidean~ Sphere}=1.
$$
 Geometrically, the fact that the quantity ~(2.74)
is independent of  vectors~$R$ means that the
indicatrix curvature is  constant. Therefore, we have arrived at

\ses
\ses

{\large   Theorem} 2.3. \it The Finsleroid Indicatrix $ \cI^{PD}_g
$ is a space of constant  positive curvature\rm.

\ses
\ses

Also, on comparing between
the result (2.74) and  Eqs. (2.22)-(2.23), we obtain

\ses
\ses

{\large  Theorem} 2.4.
 \it The Finsleroid curvature  relates to
the discriminant of the input characteristic quadratic form~\rm
(2.22) \it simply as \be \cR_{Finsleroid ~Indicatrix
}=-\fr14D_{\{B\}}. \ee \rm

\ses
\ses


Points of the indicatrix can be represented by means of the \it
unit vectors \rm $l=\{l^p\}$: \be l^p = \fr{R^p}{K(g;R)}, \ee so
that \be K(g;l)\equiv 1. \ee The vectors can conveniently be
parameterized as follows: \be l^a=n^a\fr{\sin f}h
\exp\lf(\frac12G(f-\fr{\pi}2)\rg),
 \qquad
l^N=(\cos f -\fr12G\sin f) \exp\lf(\frac12G(f-\fr{\pi}2)\rg), \ee
where \be f\in[0,\pi] \ee and $n^a$ are the components that are taken to fulfill
\be r_{ab}n^an^b=1; \ee also, \be J(g;l)=
\exp\lf(-\frac12G(f-\fr{\pi}2)\rg) \ee (cf. (2.31)). The reader is
advised to verify that
$$
A(g;l)= \fr1{J(g;l)} \cos f
$$
and
$$
\fr{hq}{A(g;l)}
=\tan f.
$$
Therefore, it is appropriate to take
 \be
  f=\arctan \fr {hq}A,
  \ee
in which case from (2.37) it follows that
$$
\Phi(g;l)=
\fr{\pi}2-f.
$$
At the same time, for the function (2.22) we find
$$
B(g;l)=
\lf(\fr1{J(g;l)}\rg)^2
=
\exp\lf(G(f-\fr{\pi}2)\rg).
$$

This method can farther be extended for the whole space by taking
the parameterizations \be R^a=\fr{K}{hJ}n^a\sin f, \qquad
R^N=\fr{K}{J} (\cos f-\fr12G\sin f), \ee \ses which entails \be
\D{R^p}K=\fr1KR^p, \ee \ses \be \D{R^a}f=\fr K{hJ}n^a (\cos
f+\fr12G\sin f), \qquad \D{R^N}f=-\fr K{h^2J}\sin f, \ee \ses \ses
\be \DD{R^a}f=\fr K{hJ}n^a \lf(G\cos f-(1-\fr14G^2)\sin f\rg), \ee
and \be \DD{R^N}f=-\fr K{h^2J} (\cos f+\fr12G\sin f). \ee

Last, we write down the explicit components of the relevant Finsleroid
Cartan tensor
\ses\\
$$
C_{pqr} := \fr12\D{g_{pq}}{R^r}:
$$
\ses
\ses
$$
R^NC_{NNN}=gw^3V^2Q^{-3}, \quad\quad R^NC_{aNN}=-gww_aV^2Q^{-3},
$$
\ses
$$
R^NC_{abN}=\fr12gwV^2Q^{-2}r_{ab}+\fr12g(1-gw-w^2)w_aw_bw^{-1}V^2Q^{-3},
$$
\ses
$$
R^NC_{abc}= -\fr12gV^2Q^{-2}w^{-1}(r_{ab}w_c+r_{ac}w_b+r_{bc}w_a)
 +gw_aw_bw_cw^{-3}\lf(\fr12Q+gw+w^2\rg)V^2Q^{-3};
$$
\ses\\
and
\ses\\
$$
R^N\3CNNN=gw^3/Q^2, \quad\quad\quad R^N\3CaNN=-gww_a/Q^2,
$$
\ses
$$
R^N\3CNaN=-gw(1+gw)w^a/Q^2,
$$
\ses
$$
R^N\3CaNb=\fr12gwr_{ab}/Q+\fr12g(1-gw-w^2)w_aw_b/wQ^2,
$$
\ses
$$
R^N\3CNab=\fr12gw\de_b^a/Q+\fr12g(1+gw-w^2)w^aw_b/wQ^2,
$$
\ses
$$
R^N\3Cabc=  -\fr12g
\lf(\de_a^bw_c+\de_c^bw_a+(1+gw)r_{ac}w^b
\rg)
/wQ
 +\fr12g(gwQ+Q+2w^2)w_aw^bw_c/w^3Q^2.
$$
\ses\\
The components have been calculated by the help of the formulae
(2.51)-(2.54).

\ses

The use of the contractions
$$
R^N\3Cabcr^{ac}=-g\fr{w^b}w\fr{1+gw}Q\lf(\fr{N-2}2+\fr1Q\rg)
$$
and
\ses\\
$$
R^N\3Cabcw^aw^c=-g\fr{w}{Q^2}(1+gw)w^b
$$
is convenient in many calculations.

Also
\ses\\
$$
R^NC_N=\fr N2gwQ^{-1}, \quad\quad R^NC_a=-\fr N2g(w_a/w)Q^{-1},
$$
\ses
$$
R^NC^N=\fr N2gw/V^2, \quad\quad R^NC^a=-\fr N2gw^a(1+gw)/wV^2,
$$
\ses
$$
C^N=\fr N2gwR^NK^{-2}, \qquad C^a=-\fr N2gw^a(1+gw)w^{-1}R^NK^{-2},
$$
\ses
\ses
\ses
$$
C_pC^p=\fr{N^2}{4K^2}g^2.
$$

\ses
\ses

The respective {\it $\cE_g^{PD}$-geodesic equation} reads
\be
\fr{d^2R^p}{ds^2}+\3Cqpr(g;R)\fr{dR^q}{ds}\fr{dR^r}{ds}=0,
\ee
where $s$ is the arc-length parameter defined by
\be
ds=\sqrt{
g_{pq}(g;R)dR^pdR^q}.
\ee

\ses
\ses

\setcounter{sctn}{3} \setcounter{equation}{0}

 2.3. {\it Associated Finsleroid Hamiltonian function}. Considering the co-vector space
$\hat V_N$ dual to the vector space $V_N$ used in the preceding
Section 2.2, and denoting by $\hat R$ the respective co-vectors,
so that $\hat R\in \hat V_N$, we may introduce the co-counterparts
of the formulas (2.1)-(2.11), obtaining the topological product
\be \hat V_N=  \hat V_{N-1}\times\hat e_N \ee and the separation
\be \hat R=\{   \hat \bR,R_N\}, \qquad R_N\in \hat e_N \quad {\rm
and} \quad \hat \bR\in \hat V_{N-1}. \ee Then we put \be R_N=\hat
Z, \ee \ses \be \hat R=\{\hat \bR,\hat Z\}, \ee and introduce a
metric \be \hat q=\hat q(\hat \bR) \ee over the
$(N-1)$-dimensional co-vector space $\hat V_{N-1}$.

With respect to a coordinate basis
$\{\hat e_a\}$
dual to
$\{e_a\}$, we obtain in $\hat V_{N-1}$ the coordinate representations
\be
\hat \bR=\{R_a\}=\{R_1,\dots,R_{N-1}\}
\ee
and
\be
\hat R=\{R_p\}=\{   R_a,R_N\}\equiv\{      R_a,\hat Z\}
\ee
together with
\be
\hat q(\hat \bR)=\sqrt{r^{ab}R_aR_b},
\ee
where $r^{ab}$ are the
contravariant
components of a symmetric positive-definte tensor
defined over $\hat V_{N-1}$;
the tensor is determined by the reciprocity $r_{ab}r^{bc}=\de_a^c$. The variables
\be
p_a=R_a/\hat Z, \qquad p^a=r^{ab}p_b,
\qquad p= \hat q/\hat Z,
\ee
where
\be
p\in(-\iy,\iy),
\ee
are convenient to apply whenever $\hat Z\ne0$.
The co-version
\be
r^{pq}=\{r^{NN}=1,~r^{Na}=0,~r^{ab}\}
\ee
of the input
metric tensor (2.11) is
meaningful over the space $\hat V_N$.
The parameter $g$ introduced in Eqs. (2.12)-(2.13), as well as the explicated
formulae (2.14)-(2.21),
are applicable in the co-approach, too.


The
{\it  characteristic
quadratic co-form}
\be
\hat B(g;R) :={\hat Z}^2-g\hat q\hat Z+{\hat q}^2
=
\fr12\Bigl[(\hat Z+g^+\hat q)^2+(\hat Z+g^-\hat q)^2\Bigr]>0
\ee
is of the negative discriminant:
\be
D_{\{\hat B\}}=-4h^2<0,
\ee
(cf. (2.22) and (2.23)).
Whenever $\hat Z\ne0$, we can use the quadratic form
\be
\hat Q(g;p) := \hat B/(\hat Z)^2,
\ee
obtaining
\be
\hat Q(g;p)=1-gp+p^2>0,
\ee
and  the function
\be
\hat E(g;p)
:=
1-\fr12gp.
\ee
Similarly to (2.26)-(2.29), we get
\be
{\hat E}^2+h^2p^2=Q,
\ee
\ses
\be
\hat B|_{_{g=0}}=r^{pq}R_pR_q,
\ee
and
\be
\hat Q|_{_{g=0}}=1+p^2.
\ee

\ses

This enables us to introduce  the {\it  Finsleroid Hamiltonian function}
(FHF for short)
\be
H(g;\hat R)=
\sqrt{\hat B(g;\hat R)}\,\hat J(g;\hat R),
\ee
where
\be
\hat J(g;\hat R)=\e^{-\frac12G\hat\Phi(g;\hat R)}
\ee
\medskip
\medskip
and
\be
\hat\Phi(g;\hat R)=
\fr{\pi}2-\arctan \fr G2-\arctan\Bigl(\fr{\hat q}{h\hat Z}-\fr G2\Bigr),
\qquad {\rm if} \quad \hat Z\ge 0,
\ee
\medskip
\be
\hat\Phi(g;\hat R)=
-\fr{\pi}2-\arctan \fr G2-\arctan\Bigl(\fr{\hat q}{h\hat Z}-\fr G2\Bigr),
\qquad {\rm if}  \quad \hat Z\le 0,
\ee
or in other forms,
\be
\hat\Phi(g;\hat R)=
\fr{\pi}2-\arctan \fr G2-\arctan\Bigl(\fr{\hat L(g;\hat R)}{h\hat Z}\Bigr),
\qquad {\rm if}  \quad \hat Z\ge 0,
\ee
\medskip
\be \hat\Phi(g;\hat R)= -\fr{\pi}2-\arctan \fr
G2-\arctan\Bigl(\fr{\hat L(g;\hat R)}{h\hat Z}\Bigr), \qquad {\rm
if}  \quad \hat Z\le 0, \ee where \be \hat L(g;\hat R)=\hat
q-\frac g2\hat Z; \ee and \be \hat\Phi(g;\hat R)= \arctan{\fr{\hat
A(g;\hat R)}{h\hat q}                  } \ee with \be \hat
A(g;\hat R)=\hat Z-\fr g2\hat q. \ee
\medskip
The respective range is such that \be -\fr{\pi}2 \le\hat \Phi\le
\fr{\pi}2 \ee
\medskip
and
\be
\hat\Phi=
\fr{\pi}2
,
\quad {\rm if} \quad \hat q=0 \quad {\rm and} \quad \hat Z>0;
\qquad
\Phi=
-\fr{\pi}2
,\quad {\rm if} \quad \hat q=0 \quad {\rm and} \quad \hat Z<0.
\ee
\medskip
Also,
\be
\hat\Phi|_{_{\hat Z=0}}=
-\arctan \fr G2.
\ee

Under these conditions, we arrive at the {\it
$\hat\cE^{PD}_g$-space}:
\be
\hat\cE^{PD}_g :=
\{\hat V_N=\hat V_{N-1}\times\hat e_N;\,\hat R\in \hat V_N;\,H(g;\hat R);\,g\}.
\ee


The function (3.21) is intimately connected with the function (2.31), namely
\be
\hat J(g;\hat R)=
\fr1{J(g;R(g;\hat R))}.
\ee

Treating the right-hand part of the definition (3.20)
as a function
$\breve H$
of  the arguments
$\{g;\hat q,\hat Z\}$,
such that
\be
\breve H(g;\hat q,\hat Z)
=
H(g;\hat R),
\ee
the co-counterparts of Eqs. (2.46)-(2.48) hold true:
\be
\breve H(g;\hat q,-\hat Z)\ne \breve H(g;\hat q,\hat Z),
\qquad unless \quad g=0;
\ee
the
\it
$g\hat Z$-parity\rm:
\be
\breve H(-g;\hat q,-\hat Z)=\breve H(g;\hat q,\hat Z);
\ee
and
\be
H(g;R)\stackrel{R_a\leftrightarrow -R_a}{\Longleftrightarrow} H(g;R).
\ee


We may also represent (3.20) in the form
\be
H(g;\hat R)=|\hat Z|W(g;p)
\ee
with the {\it  generating co-function}
\be
W(g;p)=\sqrt{\hat Q(g;p)}\,\hat j(g;p).
\ee
We have
$$
\hat j(g;p)
=
\hat J(g;1,p).
$$
Differentiating yields
\be
W'=pW/\hat Q,
\qquad
W''=W/\hat Q^2,
\ee
\ses
\be
(W^2/\hat Q)'=gW^2/\hat Q^2,  \qquad (W^2/\hat Q^2)'=2(g-p)W^2/\hat Q^3,
\ee
\ses
\be
\hat j'=\fr12g\hat j/\hat Q,
\ee
and
\be
\fr12(W^2)'=pW^2/\hat Q,
\qquad
\fr12(W^2)''=(\hat Q+gp)W^2/\hat Q^2,
\ee
together with
\be
\fr14(W^2)'''=gW^2/\hat Q^3
\ee
and
\be
\hat \Phi'=-h/\hat Q,
\ee
where the prime ($'$) denotes the differentiation with respect to~$p$.

Also,
\be
(\hat A(g;\hat R))^2+h^2{\hat q}^2=\hat B(g;\hat R)
\ee
and
\be
(\hat L(g;\hat R))^2+h^2{\hat Z}^2=\hat B(g;\hat R).
\ee


Similarly to (2.59)-(2.63),
the subsequent simple calculations yield the relations
\ses\\
$$
R^p=\fr12\D{H^2(g;\hat R)}{R_p}:
$$
\ses\ses
\be
R^N=(\hat Z-g\hat q)
\fr{H^2}{\hat B},
\qquad\quad
R^a=
r^{ab}R_b
\fr{H^2}{\hat B};
\ee
\ses
\ses
and also
\ses\\
$$
g^{pq}(g;\hat R)
=\fr12\,
\fr{\prtl^2H^2(g;\hat R)}{\prtl R_p\prtl R_q}
=\fr{\prtl R^p(g;\hat R)}{\prtl R_q}:
$$
\ses
\ses
\ses
\ses
\ses
\ses
\be
 g^{NN}(g;\hat R)=[(\hat Z-g\hat q)^2+{\hat q}^2]
\fr{H^2}{{\hat B}^2},
\qquad
g^{Na}(g;\hat R)=-g\hat q
r^{ab}R_b
\fr{H^2}{{\hat B}^2},
\ee
\ses
\ses
\ses
\ses
\be
g^{ab}(g;\hat R)=
\fr{H^2}{\hat B}
r^{ab}+g\fr{
r^{ad}R_d
r^{be}R_e
\hat Z
}{\hat q}
\fr{H^2}{{\hat B}^2}.
\ee
\ses\\
The reciprocal tensor components are
\ses\\
\ses\\
\be
g_{NN}(g;\hat R)=({\hat Z}^2+{\hat q}^2)
\fr1{H^2},
\qquad
g_{Na}(g;\hat R)=g\hat qR_a
\fr1{H^2},
\ee
\ses
\be
g_{ab}(g;\hat R)=
\fr{\hat B}{H^2}
r_{ab}-g(\hat Z-g\hat q)\fr{R_aR_b}{\hat q}
\fr1{H^2}.
\ee


To arrive from FMF (2.30)-(2.33)
at FHF (3.20)-(3.23), it is easy to note that Eq.~(2.59)
entails the equality
$$
p_a=\fr{w_a}{1+gw}
$$
which inverse is
$$
w_a=\fr{p_a}{1-gp}.
$$
Thus we find
$$
p=\fr w{1+gw}, \quad w=\fr p{1-gp}, \quad 1+gw=\fr1{1-gp},
$$
which entails the relations
$$
E(g;w)=\fr{\hat E(g;p)}{1-gp},
$$
\ses
$$
\fr{E(g;w)}{\sqrt{Q(g;w)}}=\fr{\hat E(g;p)}{\sqrt{\hat Q(g;p)}},
$$
and
\be
Q(g;w)=\fr{\hat Q(g;p)}{(1-gp)^2}=(1+gw)^2\hat Q(g;p)
\ee
together with
\be
\fr{1+gw}{Q(g;w)}=
\fr1{\sqrt{Q(g;w)\hat Q(g;p)}}.
\ee
Starting now
with (2.59), we obtain
$$
\hat Z=(Z+gq)\fr{K^2}{B}=
(Z+gq)\fr{K^2}{Z^2Q}=
\fr{1+gw}{Q}\fr{K^2}{Z}=
\fr1{\sqrt{Q\hat Q}}
\fr{K^2}{Z};
$$
whence
\be
\hat Z
=
\fr1{\sqrt{Q\hat Q}}
\fr{KH}{Z},
\ee
where the fundamental definition
$$
H(g;\hat R)=K(g;R)$$
 for the FHF has been used.
Taking into account the representations $K=ZV$ and $H=\hat ZW$
(see (2.49) and (3.38)), from (3.55) we obtain the identity \be
V^2W^2=Q\hat Q \ee which, on using $V=\sqrt{Q}j$ (see (2.50)),
just yields the FHF representation (3.20)-(3.23) together with the
equality (3.33).


Thus we have proved

\ses
\ses

{\large Theorem} 2.5.
\it The representations
\rm(3.20)-(3.23) \it associate the
required FHF
to the basic FMF given by Eqs. \rm(2.30)-(2.33).

\ses
\ses

{\large Definition}. Given the  FHF (3.20), the body
 \be
\hat\cF^{PD}_g := \{\hat R\in \hat V_N: H(g;\hat R)\le1\} \ee is
called the \it  co-Finsleroid\rm.

\ses
\ses

{\large Definition}. The respective  figuratrix defined by the
equation \be \hat\cI^{PD}_g := \{\hat R\in \hat V_N: H(g;\hat
R)=1\} \ee is called the \it co-Finsleroid Indicatrix\rm.

\ses \ses

 We remain it to the reader to verify that Theorems
2.3-2.4 proven in the preceding Section 2.2 can well be
re-formulated in the co-approach:

\ses
\ses

{\large Theorem} 2.6. \it The co-Finsleroid Indicatrix
$\hat\cI^{PD}_g$ is a constant-curvature space with the positive
curvature value~\rm(2.74): \be \cR_{co-Finsleroid~Indicatrix}=
\cR_{Finsleroid~Indicatrix} \ee \ses \it and \rm \be
\cR_{co-Finsleroid~Indicatrix}=h^2=1-\fr14g^2, \qquad\quad 0 <
\cR_{co-Finsleroid~Indicatrix} \le 1. \ee \it The formula \be
\cR_{co-Finsleroid~Indicatrix}=-\fr14D_{\{\hat B\}} \ee is valid.
\rm

\ses
\ses

We are led also to

\ses
\ses

{\large Theorem} 2.7.
\it The symmetry
\be
K(g; R)\stackrel{\left\{
\ba
{rcl}
g\longleftrightarrow -g,\\
\quad  R\longleftrightarrow \hat R\\
\ea
\right\}}
{\Longleftrightarrow}H(g;\hat R)
\ee
holds fine.\rm

\ses
\ses

{\large Note}. It is useful to verify that the insertion of the relations
(2.59) in the right-hand parts
of the representations (3.49)-(3.52) leads directly to the FMT components
given by the list (2.60)-(2.63).

\ses
\ses

\nin\setcounter{sctn}{4}\setcounter{equation}{0}

2.4. {\it Shape of Finsleroid and co-Finsleroid.} The Finsleroid
is not ``uniform" in all directions and, therefore, does not
permit general rotations. In terms of the function (2.45), the
Finsleroid equation (see the definition (2.72)) reads
 \be
  \breve K(g;q,Z)=1.
 \ee
 From (2.30)-(2.35) it follows directly that the value
 \be
  q^*
  := q_{\big|_{Z=0}}
  \ee
  of the
quantity $q$ over the Finsleroid is given by
 \be
q^*(g)=\exp\Bigl(-\fr G2\arctan{\fr G2}\Bigr);
 \ee
 with the
definitions
\be
Z_1(g)= \lf.Z\rg|_{q=0}, \quad{\rm when} \quad
Z<0,
\ee
 and
  \be
Z_2(g)  =\lf.Z\rg|_{q=0}, \quad {\rm when}\quad
Z>0,
 \ee
  we obtain
  \be
  Z_1(g)=-e^{G\pi/4} \qquad {\rm and } \qquad
Z_2(g)=e^{-G\pi/4}.
 \ee
  Thus at any given value $g$ we obtain the
simple and explicit value for the altitude of the Finsleroid:

\ses \ses

{\large Theorem} 2.8. \it We have \rm
$$
{\rm The~Altitude~of~Finsleroid}~=~Z_2(g)-Z_1(g)=2\cosh\fr{G\pi}4.
$$

\ses
\ses

 The equation (4.1) cannot be resolved  to find  the function
\be
Z=Z(g;q)
\ee
 in an explicit form, because of a rather high complexity
of the right-hand parts of Eqs. (2.30)-(2.35). Nevertheless,
differentiating the identity
 \be
 \breve K(g;q,Z(g;q))=1
  \ee
   (see
(4.1))  yields, on using (2.51), the simple results for the first
derivatives:
 \be \D{Z(g;q)}{q}=-\fr q{Z+gq}
  \ee
   and \be
\DD{Z(g;q)}{q}=-\fr{B(g;R)}{(Z+gq)^3}. \ee
\bigskip
We also get
\be
\lf.\fr{dZ(g;q)}{dq}\rg|_{q=0}=0\qquad {\rm and}  \qquad
\fr{dZ(g;q)}{dq}\mathop{\Longrightarrow}\limits_{Z\to+0}-\fr1g.
\ee

Inversely, for the function
\be
q=q(g;Z)
\ee
obeying (4.1) we obtain
\ses
\be
\D{q}{Z}=-g-\fr Zq
\ee
\ses\\
and \ses \be \DD{q}{Z}=-\fr{B(g;R)}{q^3}<0. \ee
\ses\\
We have \ses
\be \D{q(g;Z)}{Z}>0, \quad {\rm if} \quad  Z<-gq;
 \qquad \quad
\D{q(g;Z)}{Z}<0, \quad {\rm if} \quad  Z>-gq.
\ee

Also,
\ses
\be
\D{q}{Z}=0,
\quad \rm if \it\quad  Z=  Z^{**} \quad \rm with \it\quad   Z^{**}=-gq^{**}.
\ee
\ses\\
Inserting this $ Z^{**}$ in (2.30)-(2.35) yields \ses \be
 \Phi^{**}=
\fr{\pi}2+\arctan \fr G2-\arctan {\fr{g^2-2}{2gh}}, \qquad {\rm
if} \quad Z^{**}\ge 0 ~\sim ~g<0, \ee
\medskip
\be \Phi^{**}= -\fr{\pi}2+\arctan \fr G2-\arctan
{\fr{g^2-2}{2gh}}, \qquad {\rm if} \quad Z^{**}\le 0~ \sim ~g>0,
\ee
\ses\\
and
\be
 q^{**}(g)=
\e^{-\frac12G\Phi^{**}}
\ee
\ses\\
together with
\be
 Z^{**}(g)=-g
\e^{-\frac12G\Phi^{**}}.
\ee

Therefore, the following assertion can be set up for the width
 of the Finsleroid.

\ses\ses

{\large Theorem} 2.9. \it With any given $g$, \rm
$$
{\rm The~Width~of~Finsleroid}~=~2q^{**}(g) =
2\e^{-\frac12G\Phi^{**}}.
$$

\ses\ses

The formulas (4.19) and (4.20) may also be interpreted by saying
that The {\it Equatorial Section} of the Finsleroid is of the
radius \be r_{Equatorial}=q^{**}\ee and cuts the $Z$-axis at \be
Z_{Equatorial}=Z^{**}. \ee
With the parameter value $|g|$ being
increasing, the Finsleroid is stretching  in wide and altitude:
 \be
q^{**}\,\mathop{\Longrightarrow}\limits_{|g|\to 2}\,\infty \ee and
\be |Z^{**}|\,\mathop{\Longrightarrow}\limits_{|g|\to 2}\,\infty,
\ee tending in its shape to a cone:
 \be
\fr{q^{**}}{|Z^{**}|}\,\mathop{\Longrightarrow}\limits_{|g|\,\to
2} \fr12,
 \ee
such that the vertex of the Finsleroid tends to approach the
origin point ``O". From (4.23) one can infer

\ses\ses

{\large Theorem} 2.10. \it We have: \rm
$$
{\rm The~Limiting~Vertex~Angle}~=~2\arctan\fr12.
$$

\ses\ses

The above formulae, particularly the negative sign of
the second derivative (4.10), can be used directly to verify
the following

\ses\ses

{\large Theorem} 2.11. \it The Finsleroid Indicatrix $\cI_g^{PD}$
is closed, regular, and strongly convex.\rm

\ses \ses

The co-Finsleroid equation \be
 \breve H(g;\hat q,\hat Z)=1
 \ee
(cf. Eqs. (3.35) and (3.58)) can be studied in a similar way,
leading to the relations obtainable from Eqs. (4.3)-(4.20) by
means of the formal replacement $ \{g\rightarrow -g,
 R\rightarrow \hat R\}
$, owing to the fundamental symmetry (3.61).

Therefore, we can state the following:

\ses\ses

{\large Theorem} 2.12. \it The co-Finsleroid Indicatrix
$\hat\cI^{PD}_g$ is closed, regular, and also strongly convex.\rm

\ses\ses

{\large Theorem} 2.13. \it At any given parameter $g$, the
Finsleroid and
 the co-Finsleroid mirror one another
under the $g$-reflection: \rm \be
\cF^{PD}_g\stackrel{g\longleftrightarrow -g}
{\Longleftrightarrow}\hat\cF^{PD}_{-g}. \ee

\ses
\ses

\ses
\ses

All  Figures shown below have been prepared by means of a precise
use of Maple9.

In Figs.2-7  bold lines serve to draw the Finsleroids, while unit
circles  simulate the ordinary Euclidean spheres. Fig.2 may be
used as a convenient demonstration example  (the \it  trainer\rm)
for the Finsleroid by showing various structure details, including
the equatorial section and the characteristic tangents, in a
distinct way. We remain it to the reader to evaluate the angles
that are depicted in the example and find among them equal cases.

Figs.2-7 clearly support the validity of Theorem 2.11 about
regularity and convexity and make an idea of existence a
diffeomorphic spherical map (see Eq. (5.1) below) a quite
trustworthy one.

If one compares between  Fig.2 and Fig.3 between Fig.4and Fig.5
or  between   Fig.~6 and Fig.7, one observes immediately that the
change of sign of the characteristic Finslerian parameter $g$ does
{\it turn up} the figures and, therefore, verifies the fundamental
Finslerian $Z$-parity property (as given by Eq. (2.47)) in a due
visual way. In a narrow sense, Figs.2-7 show the geometry of the
generatrix for the Finsleroids, the latter being (hyper)surfaces
of revolution over the $Z$-axis.

It can be traced also how the parameter $g$ effects the shape of
Finsleroid. A positive value $g$ (deforms and) shifts the unit
sphere in the down-wise manner, respectively a negative value in
the up-wise manner.

 Fig.8 and  Fig.9 model the
important functions (4.3) and (4.20), respectively.


There is no need to picture  co-Finsleroids, for they mirror
Finsleroids  with respect to the $(R^N=0)$-plane (according to our
Theorems 2.12 and 2.13 and corresponding Eqs. (3.62) and (4.27)).
In particular, at $g=-0.4$, the co-Finsleroid looks like the
demonstration example given by Fig.1 at $g=0.4$ (with
interchanging respectively  the coordinate axes: $R^N$ with $P_N$
and $\bR$ with $\bP$).


\pgbrk

\begin{figure}[!ht]
\centering
\includegraphics[width=12cm]{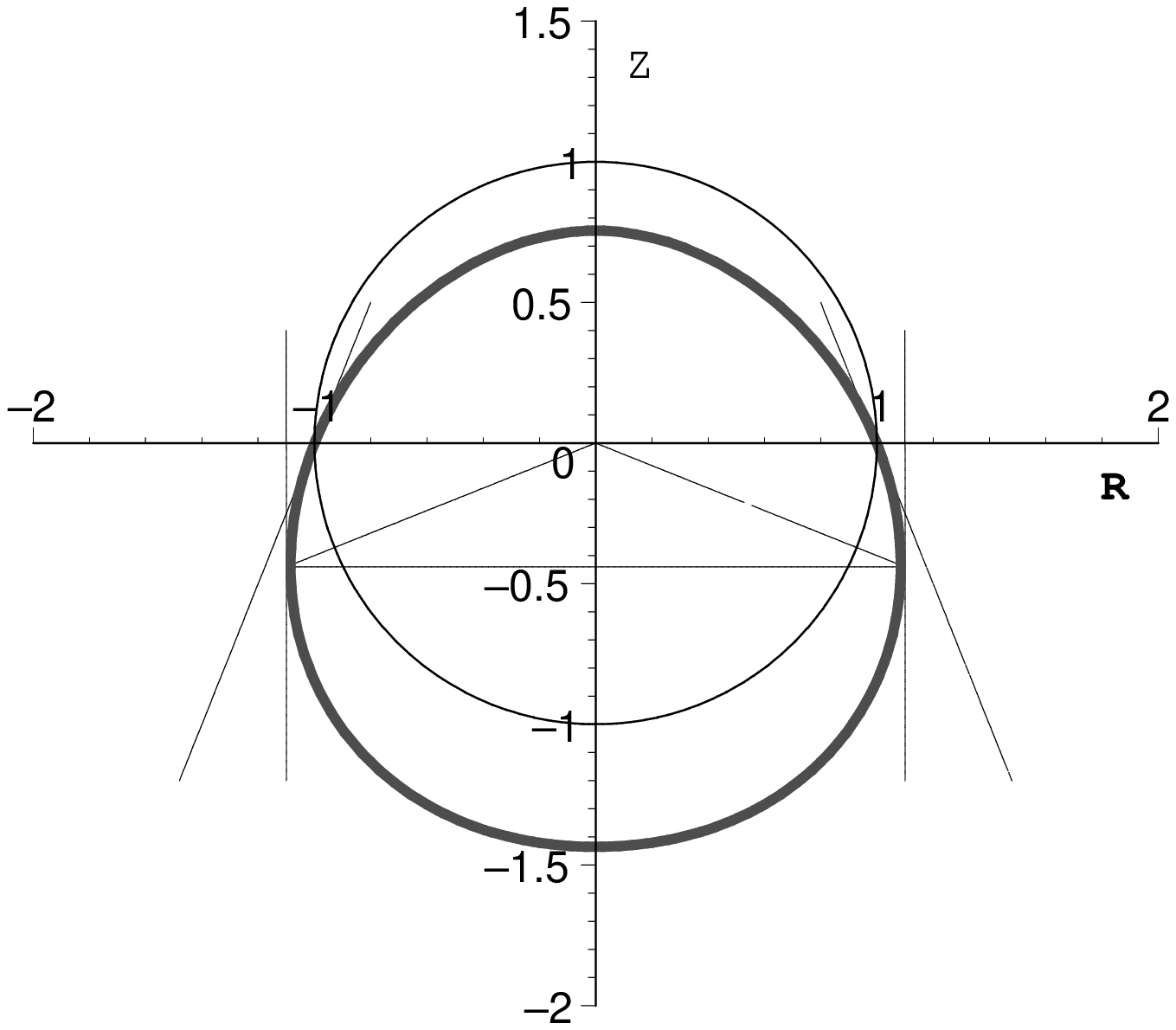}
\caption{[$g$=0.4]} \vspace{1cm}
\includegraphics[width=12cm]{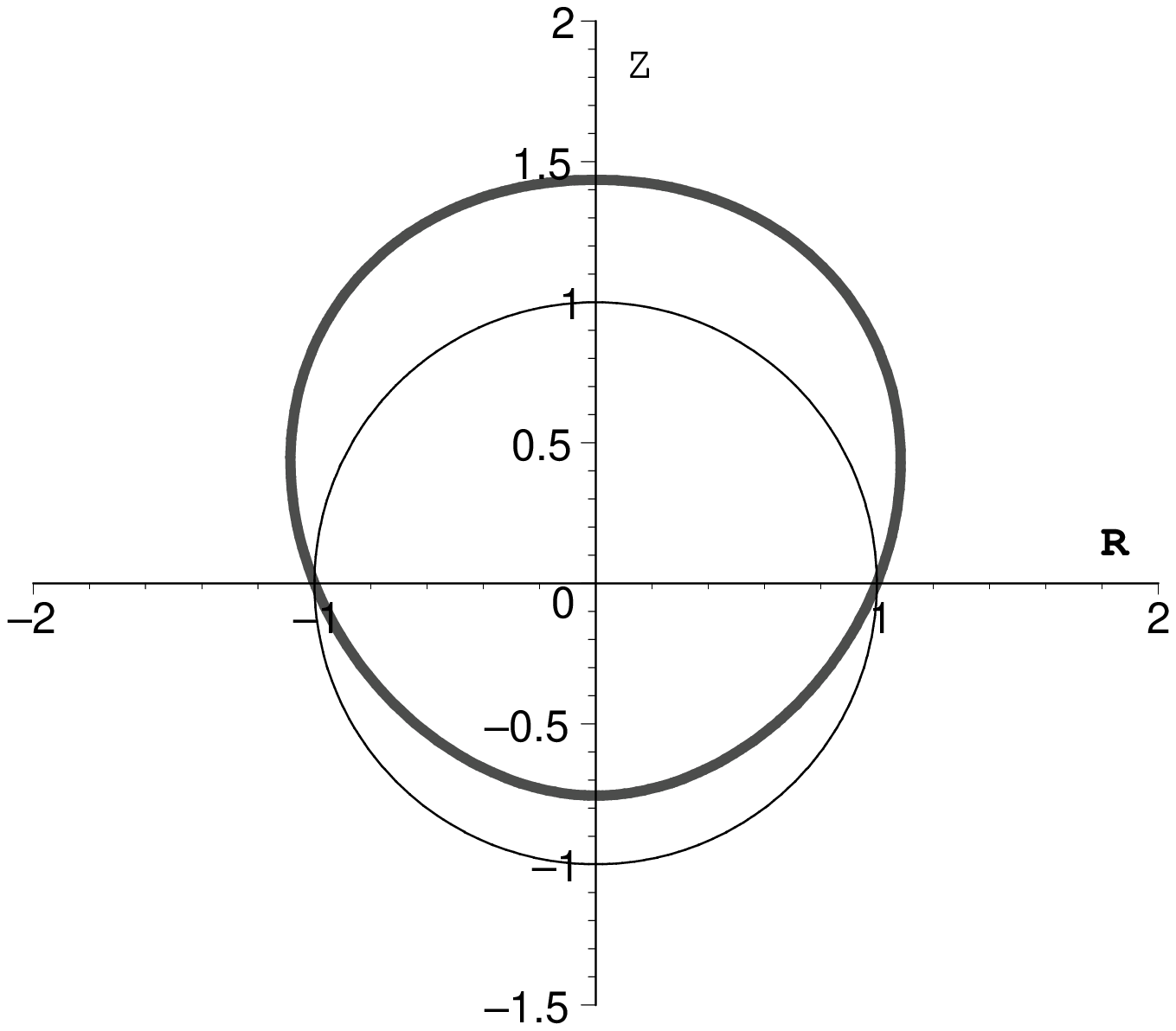}
\caption{[$g$=-0.4]}
\end{figure}

\pgbrk

\begin{figure}[!ht]
\centering
\includegraphics[width=10cm]{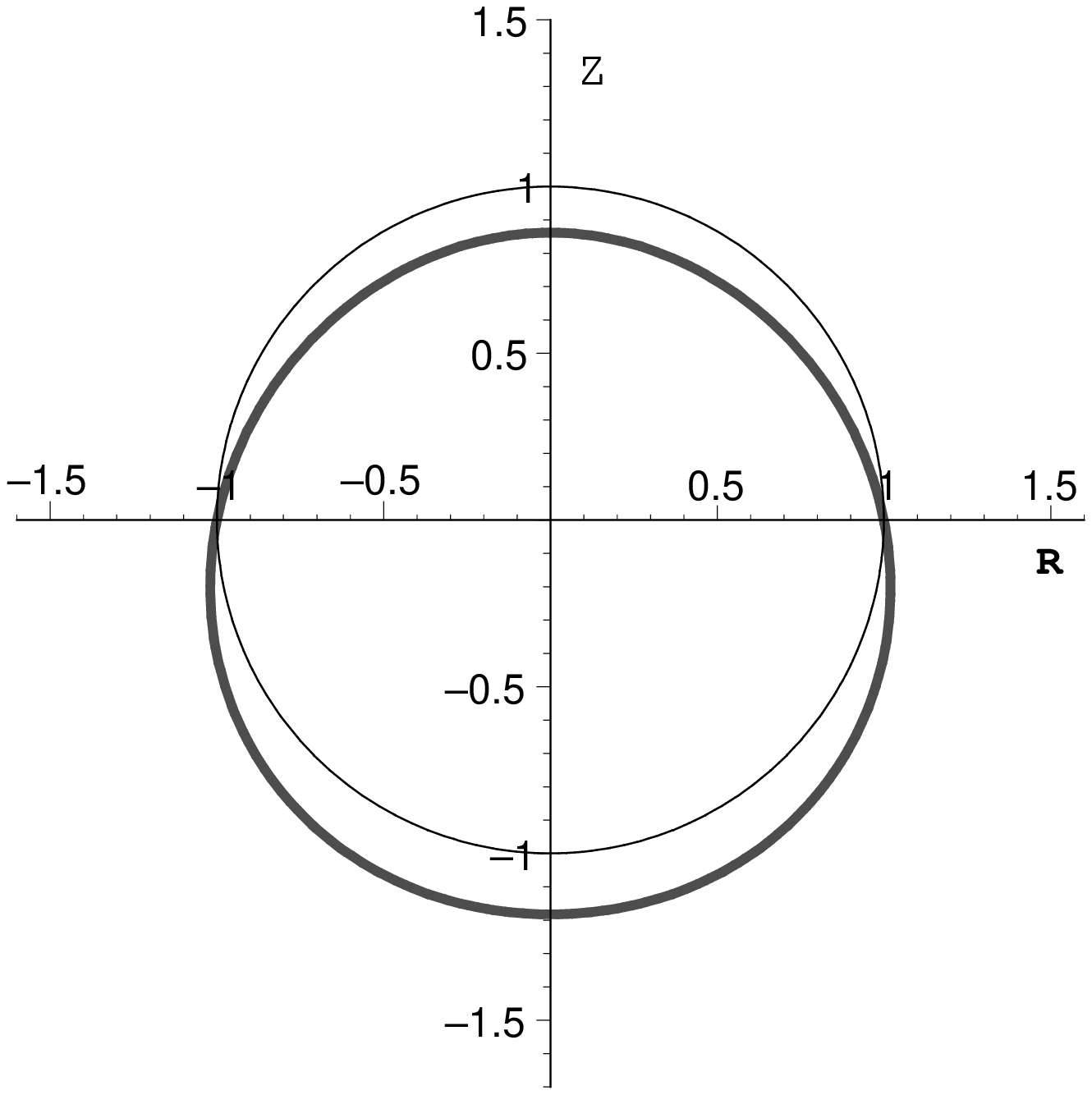}
\caption{[$g$=0.2]} \vspace{1cm}
\includegraphics[width=10cm]{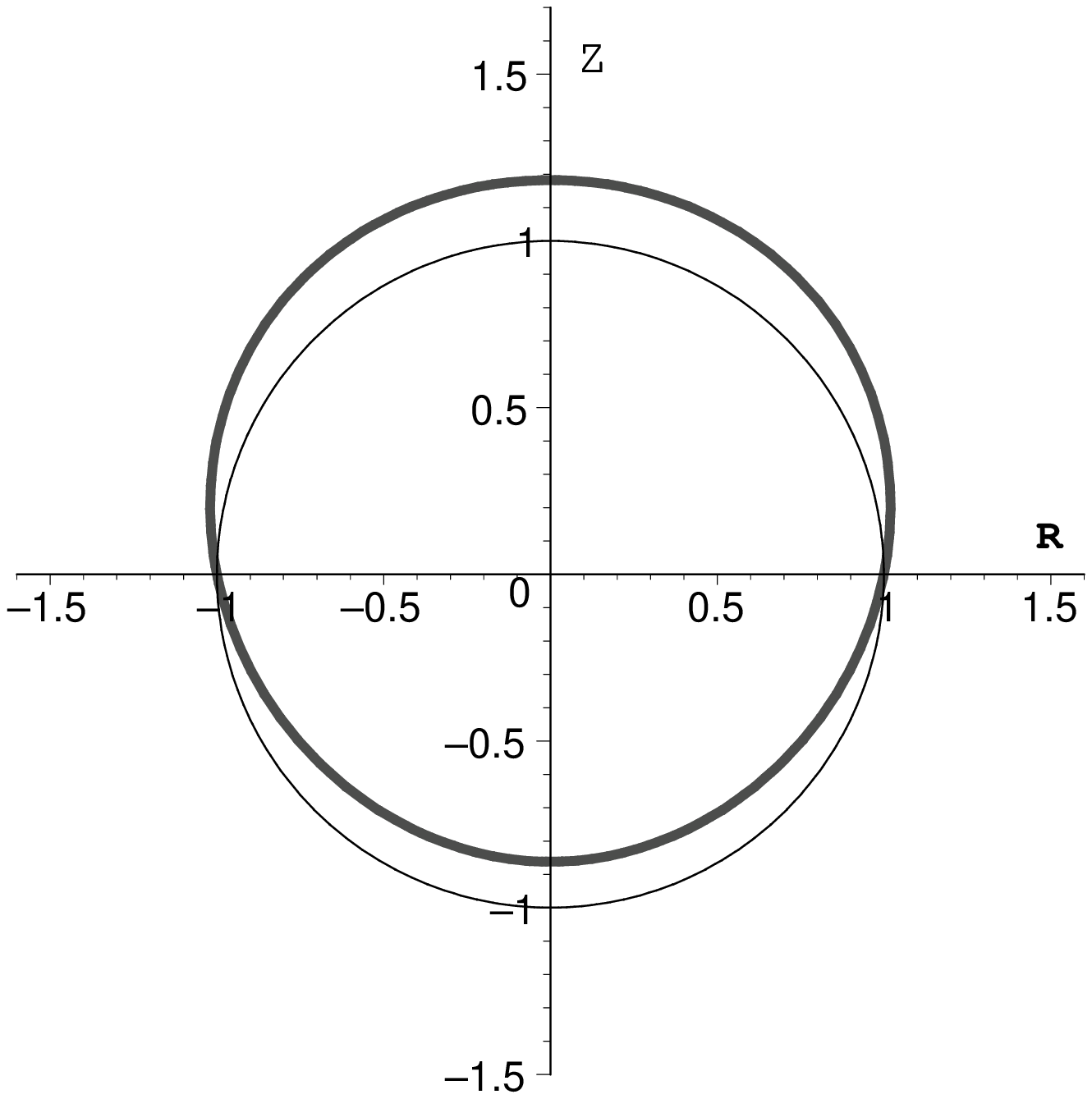}
\caption{[$g$=-0.2]}
\end{figure}

\pgbrk

\begin{figure}[!ht]
\centering
\includegraphics[width=10cm]{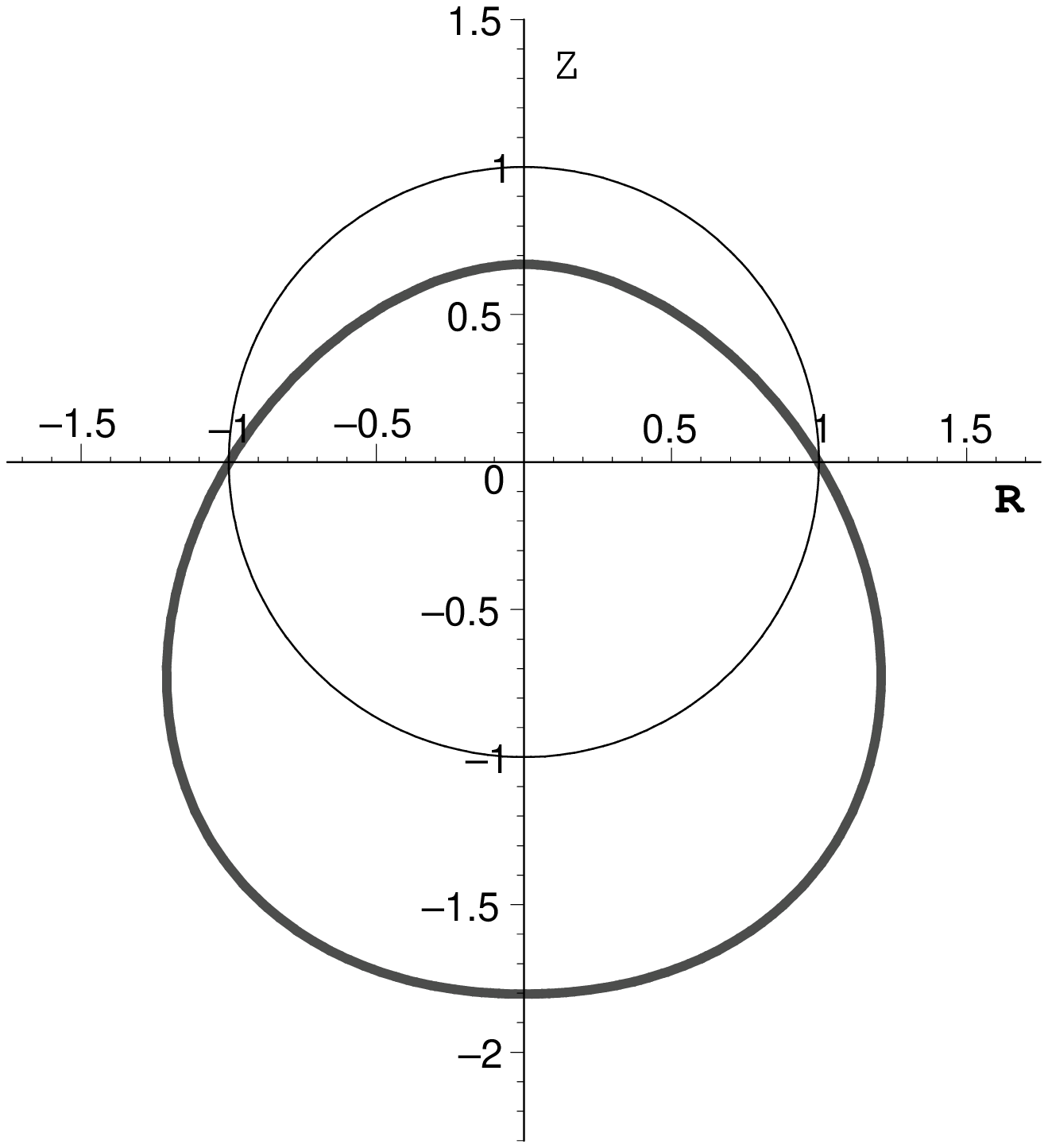}
\caption{[$g$=0.6]} \vspace{1cm}
\includegraphics[width=10cm]{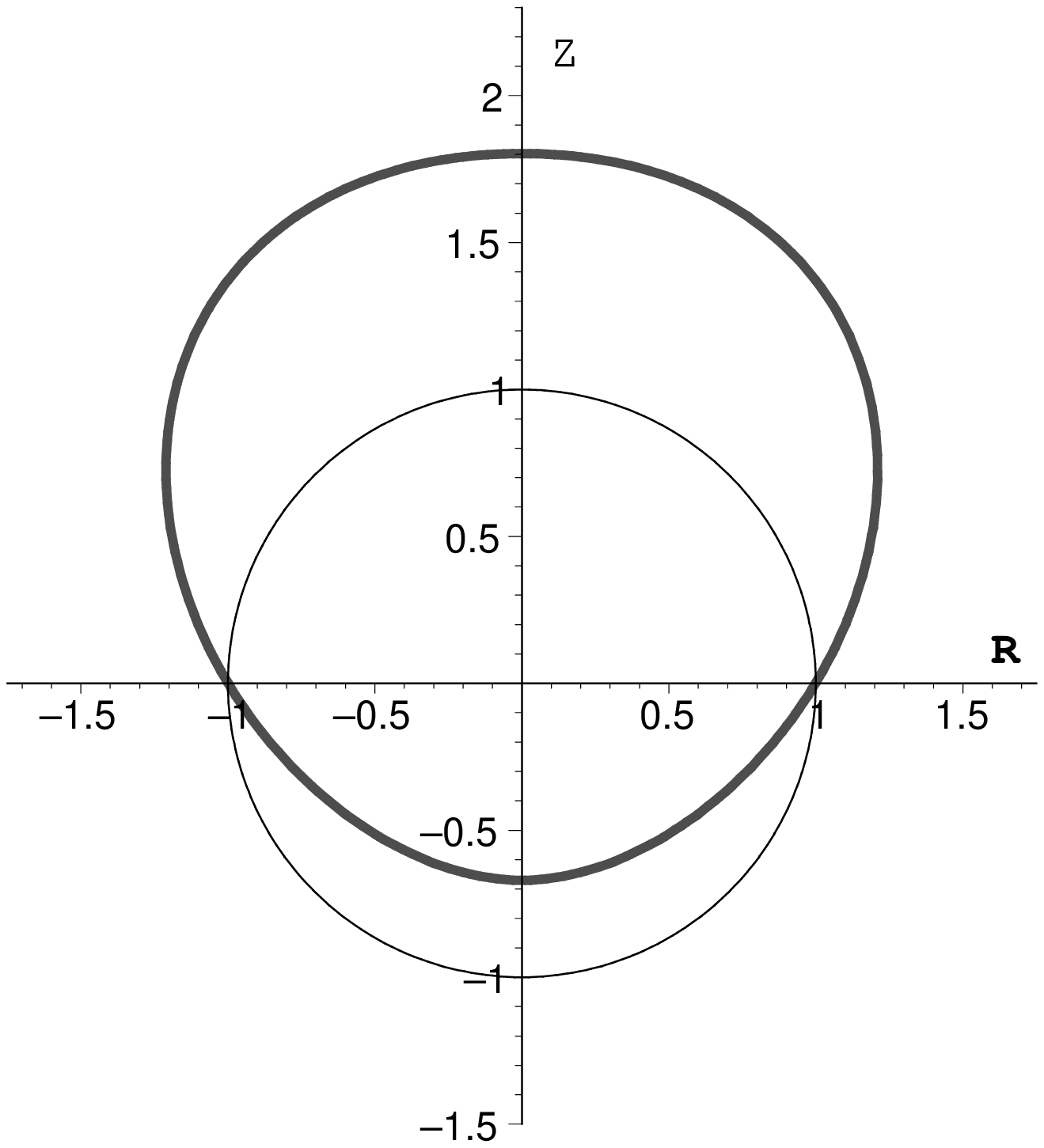}
\caption{[$g$=-0.6]}
\end{figure}

\pgbrk

\begin{figure}[!ht]
\centering
\includegraphics[width=8cm]{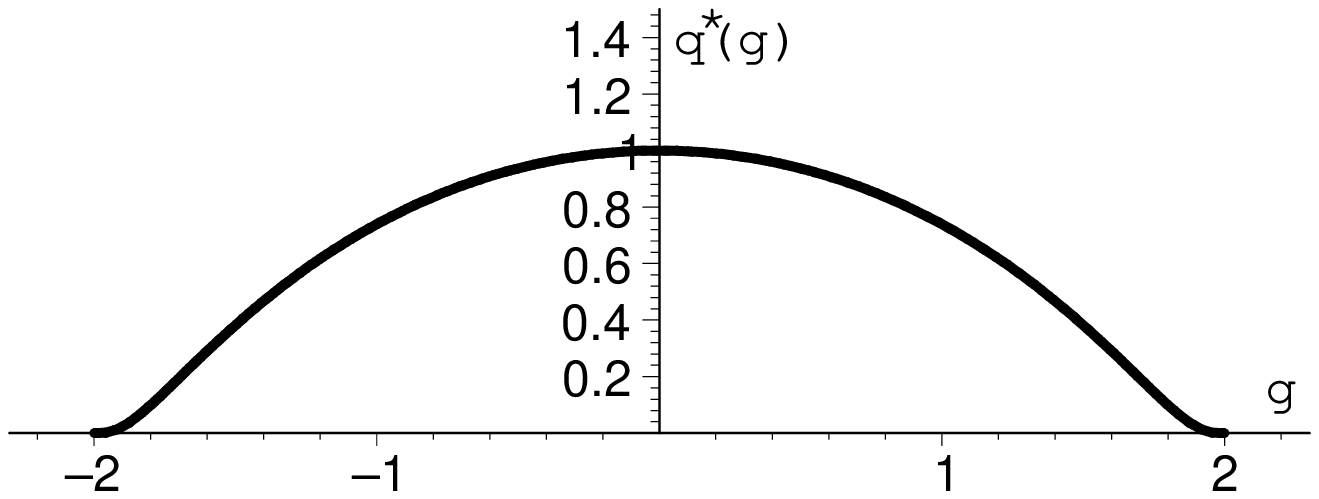}
\caption{[Eq. (4.3)]}\vspace{4cm} \centering
\includegraphics[width=8cm]{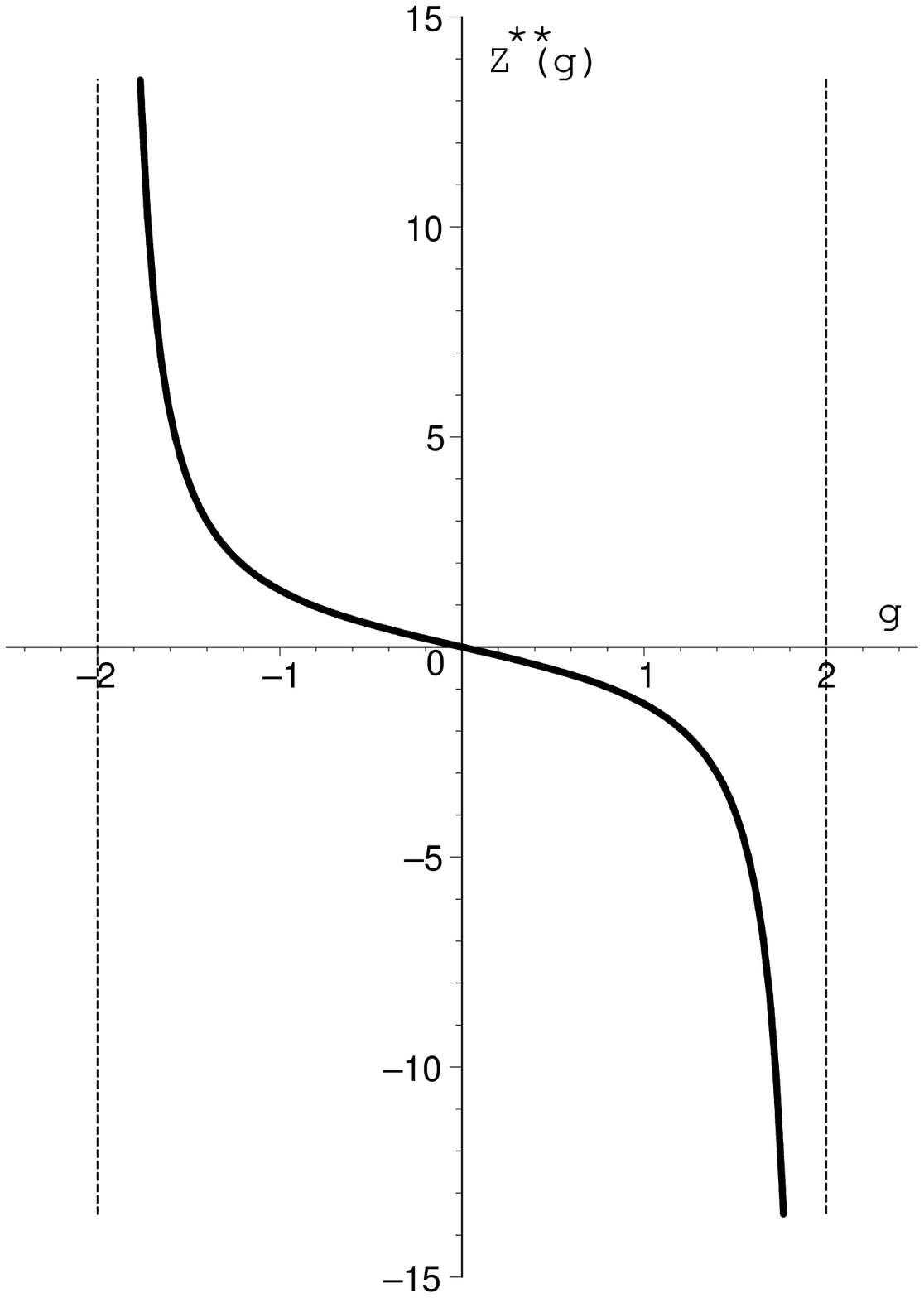}
\caption{[Eq. (4.20)]}
\end{figure}

\pgbrk










\clearpage

\setcounter{sctn}{5} \setcounter{equation}{0}

2.5. {\it Quasi-Euclidean  map of Finsleroid.} Theorem 2.11 can be
continued farther by indicating  the diffeomorphism \be \cF_g^{PD}
\stackrel{{\it i_g}}{\Longrightarrow}\cB^{PD} \ee
 of
the Finsleroid $\cF_g^{PD}\subset V_N$ to the unit ball
$\cB^{PD}\subset V_N$: \be \cB^{PD} :=\{R\in\cB^{PD}:\;
S(R)\le1\}, \ee
 where \be S(R)=\sqrt{r_{pq}R^pR^q} \equiv
\sqrt{(R^N)^2+r_{ab}R^aR^b} \ee is the input Euclidean metric
function (see (2.11)).

The diffeomorphism (5.1) can always be extended to get the
diffeomorphic map \be V_N \stackrel{\si_g}{\Longrightarrow}V_N \ee
of the whole vector space $V_N$ by means of the homogeneity: \be
\si_g\cdot (bR)=b\si_g\cdot R, \quad b>0. \ee To this end it is
sufficient to take merely \be \si_g\cdot R=||R||
i_g\cdot\Bigl(\fr{R}{||R||}\Bigr), \ee where \be ||R||=K(g;R). \ee
Eqs. (5.1)-(5.7) entail \be K(g;R)=S(\si_g\cdot R). \ee

On the other hand, the identity (2.57) suggests  taking the map
\be \bar R=\si_g\cdot R \ee by means of the components \be
 \bar R^p=\si^p(g;R)
\ee
with
\ses
\be
\si^a=
R^ah
J(g;R),
\qquad
\si^N=A(g;R)J(g;R),
\ee
\ses\\
where $J(g;R)$ and $A(g;R)$ are the functions (2.31) and (2.38).
Indeed, inserting (5.11) in (5.3) and taking into account
Eqs. (2.30) and (2.57),
we get
the  identity
\be
S(\bar R)=K(g;R)
\ee
which is tantamount to the implied relation
(5.8).

Thus we have arrived at

\ses\ses

{\large Theorem}  2.14. \it The map given explicitly by Eqs.
\rm(5.9)-(5.11) \it assigns the diffeomorphism between the
Finsleroid and the unit ball according to Eqs. \rm(5.1)-(5.8).

\ses\ses

This motivates the following

\ses\ses

{\large  Definition.}  Under these conditions, the map (5.1) is
called the { \it  quasi-Euclidean map of Finsleroid.}

\ses \ses

The inverse
\ses
\be
 R=\mu_g\cdot \bar R,
  \qquad \mu_g=(\si_g)^{-1}
\ee
\ses\\
of the transformation (5.9)-(5.11)
can be presented by the components
\be
  R^p=\mu^p(g;\bar R)
\ee
with
\be
\mu^a=\bar R^a/hk(g;\bar R),
\qquad
\mu^N=I(g;\bar R)/k(g;\bar R),
\ee
\ses\\
where
\ses
\be
k(g;\bar R)
:= J(g;\mu(g;\bar R))
\ee
\ses\\
and
\be
I(g;\bar R)
=\bar R^N-\frac12G\sqrt{r_{ab}\bar R^a\bar R^b}.
\ee
The identity
\ses
\be
\mu^p(g;\si(g;R))\equiv R^p
\ee
\ses\\
can readily be verified. Notice that
 \ses
 $$ \fr{\bar R^a}{S(\bar
R)}=\fr{hR^a}{\sqrt{B(g;R)}}, \qquad \fr{\bar R^N}{S(\bar
R)}=\fr{A(g;R)}{\sqrt{B(g;R)} },$$ \ses \be
 \fr{\sqrt{r_{ab}\bar
R^a\bar R^b}}{\bar R^N} =\fr{hq}{A(g;R)}, \qquad
w^a=\fr{R^a}{R^N}=\fr{\bar R^a}{hI(g;\bar R)},
 \ee
\ses\\
and
\be \sqrt{B}/R^N=S/I, \qquad \sqrt{Q}=S/I. \ee

The $\si_g$-image \ses \be \phi(g;\bar R)
 :=\Phi(g;R)|_{_{R=\mu(g;\bar R)}} \ee
\ses\\
of the function $\Phi$ described by Eqs. (2.32)-(2.42) is of a
clear meaning of  angle:
 \ses \ses
 \be \phi(g;\bar
R)= \arctan\fr{\bar R^N}{\sqrt{r_{ab}\bar R^a\bar R^b}}
\ee
\ses\\
(Eq. (5.19) has been used)
which ranges over \ses \ses \be
-\fr{\pi}2 \le\phi\le \fr{\pi}2. \ee
\ses\\
We have
\ses
\be
\phi=
\fr{\pi}2,
\quad {\rm if} \quad \bar R^a=0 \quad {\rm and} \quad \bar R^N>0;
\qquad
\phi=
-\fr{\pi}2
,\quad {\rm if} \quad \bar R^a=0 \quad {\rm and} \quad \bar R^N<0,
\ee
\ses\\
and also
\ses
\be
\phi|_{_{\bar R^N=0}}=0.
\ee
\ses\\
Comparing Eqs. (5.16) and (2.31) shows that
\ses
\be
k=\e^{\frac12G\phi}.
\ee

The right-hand parts in (5.11) are homogeneous
 functions
of degree~1:
\be
\si^p(g;bR)=b\si^p(g;R), \quad b>0.
\ee
Therefore, the identity
\be
\si_s^p(g;R)R^s=\bar R^p
\ee
should be valid for the derivatives
\ses\\
\be
\si_p^q(g;R)
:= \D{\si^q(g;R)}{R^p}.
\ee
The simple representations
\be
\si^N_N(g;R) =\lf(B+\fr12gqA\rg)\fr{J}{B},
\ee
\bigskip
\be
\si^N_a(g;R) =-\fr{g(ZA-B)}{2q}\fr{J r_{ab}R^b}{B},
\ee
\bigskip
\be
\si^a_N(g;R) =\fr12gq\fr{J R^ah}{B},
\ee
\bigskip
\be
\si^a_b(g;R) =\lf(B\delta^a_b-
\fr{gr_{bc}R^cR^aZ}{2q}\rg)\fr{Jh}{B},
\ee
\ses\\
and also the determinant
\ses
\be
\det(\si_p^q)=h^{N-1}J^N
\ee
\ses\\
are obtained. The relations
$$
\si_b^aR^b=JhR^a(AZ+q^2)/B, \qquad
r^{cd}\si_c^a\si_d^b=
J^2h^2\lf[r^{ab}-g(R^aR^bZ/qB)+\fr14g^2(R^aR^bZ^2/B^2)\rg]
$$
are handy in many calculations involving the coefficients $\{\si_p^q\}$.

Henceforth, to simplify notation, we shall use the substitution
\be
t^p=\bar R^p.
\ee

Again, we can note the homogeneity
\be
\mu^p(g;bt)=b\mu^p(g;t), \quad b>0,
\ee
for the functions (5.15), which entails the identity
\be
\mu_s^p(g;t)t^s=R^p
\ee
for the derivatives
\ses
\be
\mu_q^p(g;t)
:=\D{\mu^p(g;t)}{t^q}.
\ee
\ses\\
We find
\ses
\be
\mu_N^N=1/k(g;t)-\fr12g
\fr{m(t)
I(g;t)}{k(g;t)(S(t))^2},
\qquad
\mu_a^N=\fr12g\fr{r_{ac}t^cI^*(g;t)}{k(g;t)(S(t))^2},
\ee
\ses
\ses
\be
\mu_N^a=-\fr12g
\fr{m(t)
\,t^a}{hk(g;t)(S(t))^2},
\qquad
\mu_b^a=\fr1{hk(g;t)}\de_b^a+
\fr12g\fr{t^Nt^ar_{bc}t^c}
{m(t)
\,hk(g;t)(S(t))^2},
\ee
\ses\\
where
\be
m(t)=
\sqrt{r_{ab}t^at^b},
\ee
\ses
\be
I^*(g;t)=hm(t)-\fr12gt^N,
\ee
and
\ses
\be
S(t)=\sqrt{r_{rs}t^rt^s}
\equiv
\sqrt{(t^N)^2+r_{ab}t^at^b}.
\ee

\ses The  relations
$$
\D{(1/k(g;t))}{t^N}=-\fr12g\fr{m(t)}{hk(g;t)(S(t))^2},
\quad\quad
\D{(1/k(g;t))}{t^a}=\fr12g\fr{t^Nr_{ab}t^b}
{m(t)hk(g;t)(S(t))^2}
$$
\ses\\
are obtained.

Also,
\be
 R_p\mu^p_q=t_q,
\qquad
t_p\si^p_q=R_q.
\ee

The unit vectors
\be L^p
:= \fr{t^p}{S(t)},
 \qquad L_p
 := r_{pq}L^q \ee fulfill the relations \be L^q=l^p\si^q_p, \qquad
l^p=\mu^p_qL^q, \qquad l_p=\si^q_pL_q, \qquad L_p=\mu^q_pl_q, \ee
where $l^p=R^p/K(g;R)$ and $l_p=R_p/H(g;R)$ are the initial
Finslerian unit vectors.

Now we use the explicit formulae (2.62)-(2.63) and (5.30)-(5.33)
to find the transform \be n^{rs}(g;t) := \si_p^r\si_q^sg^{pq} \ee
of the FMT  under the $\cF_g^{PD}$-induced map (5.9)-(5.11), which
results in

\ses\ses

{\large Theorem} 2.15. \it One obtains the simple representation
\ses

\be
n^{rs}=h^2r^{rs}+\fr14g^2L^rL^s.
\ee
The covariant version reads
\be
n_{rs}=
\fr1{h^2}r_{rs}-\fr14G^2L_rL_s.
\ee
The determinant of this tensor is a constant:\rm
\be
\det(n_{rs})=h^{2(1-N)}\det(r_{ab}).
\ee

\ses\bigskip

Notice that
$$
L^pL_p=1, \quad n_{pq}L^q=L_p, \quad n^{pq}L_q=L^p, \quad
n_{pq}L^pL^q=1, \quad n_{pq}t^pt^q=(S(t))^2.
$$
\ses

Eq. (5.47) obviously entails
\be
g_{pq}=
n_{rs}(g;t)\si_p^r\si_q^s.
\ee






\clearpage

\setcounter{sctn}{6}
\setcounter{equation}{0}

2.6. {\it  Quasi-Euclidean metric tensor}.
Let us introduce

\ses
\ses

{\large Definition}.  The metric tensor $\{n_{pq}(g;t)\}$ of the form (5.48) is called {\it
quasi-Euclidean}.

\ses
\ses

{\large Definition}.  The \it
 quasi-Euclidean space
\be
\cQ_N :=\{V_N;n_{pq}(g;t);g\}
\ee
\rm
is an extension of the Euclidean space $\{V_N;r_{pq}\}$ to the case $g\ne0$.

\ses
\ses

The transformation (5.47) can be inverted to read
 \be
g_{pq}=\si_p^r\si_q^sn_{rs}. \ee For the angular metric tensor
(see the formula going below Eq. (2.65) in Chapter 2), from (5.46)
and (6.2) we infer \be h_{pq}=\si_p^r\si_q^sH_{rs}\fr1{h^2}, \ee
where \be H_{rs} := r_{rs}-L_rL_s \ee is the tensor showing the
orthogonality property \be L^rH_{rs}=0. \ee \ses One can readily
verify that
$$
H_{rs}= h^2(n_{rs}-L_rL_s).
$$

\ses
\ses

{\large Theorem} 2.16.  \it
The quasi-Euclidean metric tensor
is conformal to the  Euclidean metric
tensor.
\rm

\ses
\ses

Indeed, if we consider the map
\be
\bar R^p\rightarrow \tilde
R:\qquad \tilde R^p=\xi(g;\bar R)\bar R^p/h
\ee
with
 \be
 \xi(g;\bar
R)=a\lf(g;\,\fr12S^2(\bar R)\rg)
 \ee
and use the coefficients
 \be
k_q^p := \D{\tilde R^p}{\bar R^q} =(\xi\de_q^p+a'\bar R^p\bar
R_q)/h \ee
 to define the tensor
\bigskip
\be c^{pq}(g;\tilde R) := k_r^pk_s^qn^{rs}(g;\bar R), \ee we find
that \be
 c^{pq}=\xi^2r^{pq}
\ee
whenever
\be
\xi=\lf[\fr12S^2(\bar R)\rg]^{(h-1)/2}.
\ee
 The proof of Theorem 2.16
is complete.

\ses
\ses
\ses
\ses

The respective determinant representation is found to be merely
\be
 \det(c_{pq})=\xi^{-2N}\det(r_{pq}).
\ee


Let us now use the obtained
quasi-Euclidean metric tensor
$
n_{pq}(g;t)$ to construct the associated
{\it quasi-Euclidean Christoffel symbols}
$\3N prq(g;t)$. We find consecutively:
\be
n_{pq,r}
:= \D{n_{pq}}{t^r}=-\fr14G^2(H_{pr}L_q+H_{qr}L_p)/S,
\ee
and
\be
\3N prq=n^{rs}N_{psq}, \qquad
N_{prq}=\fr12(n_{pr,q}+n_{qr,p}-n_{pq,r}),
\ee
together with
\be
N_{prq}(g;t)=-\fr14G^2H_{pq}L_r/S,
\ee
which eventually yields
\be
\3N prq(g;t)=-\fr14G^2L^rH_{pq}/S.
\ee
Comparing the representation (6.16) with the identity (6.5) shows that
\be
t^p\3N prq=0, \qquad \3N pss=0,\qquad
\3 N tsr \3N ptq=0.
\ee

Also,
\be
\D{\3N prq}{t^s}-\D{\3N prs}{t^q}
 =-\fr14G^2(H_{pq}{H_s}^r-H_{ps}{H_q}^r)/S^2.
\ee Using the identities (6.17)-(6.18) in the {\it quasi-Euclidean
curvature tensor}: \be R_p{}^r{}_{qs}(g;t)
:=\D{N_p{}^r{}_q}{t^s}-\D{N_p{}^r{}_s}{t^q}+
N_p{}^w{}_qN_w{}^r{}_s-N_p{}^w{}_sN_w{}^r{}_q, \ee we arrive at
the simple result:
\bigskip
\be
R_{prqs}(g;t)     =-\fr14G^2(H_{pq}H_{rs}-H_{ps}H_{qr})/S^2.
\ee
This infers the identities
\be
L^pR_{pqrs}=L^qR_{pqrs}=L^rR_{pqrs}=L^sR_{pqrs}=0.
\ee

\ses
\ses

{\large Note}. Because of the transformation rules (5.11) and (5.47),
the representation (6.20) is tantamount to
Eqs. (2.70)-(2.71).
Therefore {\it we have got another rigorous proof
of Theorem } 2.3, {\it and of
Eq.} (2.74), {\it concerning the Finsleroid curvature}.

\ses
\ses


The {\it  quasi-Euclidean  orthonormal frames}
defined by the representations:
\be
n^{pq}=\suml_{P=1}^{N}m_P^pm_P^q
\ee
and
\be
n_{pq}=\suml_{P=1}^{N}f^P_pf^P_q
\ee
prove to be taken as
\ses\\
\be
f_q^P(g;t)=\fr1hh_q^P+\fr{h-1}hL_qL^P
\ee
and
\be
m_P^q(g;t)=hh_P^q+(1-h)L_PL^q
\ee
with
\be
L^P=h_q^PL^q, \qquad L_P=h_P^qL_q.
\ee
\ses\\
Here, $h^P_q$ and $h_P^q$
are the orthonormal frames for the input metric tensor (2.11):
\be
r^{pq}=\suml_{P=1}^{N}h_P^ph_P^q
\ee
and
\be
r_{pq}=\suml_{P=1}^{N}h^P_ph^P_q.
\ee
We have
\be
f_q^P(g;t)t^q=t^P.
\ee
\ses

The {\it  associated quasi-Euclidean Ricci rotation coefficients}:
\be
{R^{PQ}}_p(g;t)=\lf(\prtl_pf_q^Q-\3Nprqf_r^Q\rg)m_T^q\de^{TP}
\ee
are found in the simple
explicit form
\be
{R^{PQ}}_p(g;t)=(h-1)(L^Pf_p^Q-L^Qf_p^P)/S(t).
\ee


The structure of the right-hand part of the representation (6.20)
is such that  the quasi-Euclidean metric tensor $n_{pq}$ relates to a
space of  constant curvature only in the Euclidean case $g=0$ proper.
At the same time, treating the $r$-radius sphere
\be
\cS_r :=\{R\in\cS_r:\, S(t)=r\}
\ee
(the metric (5.3) has been used)
considered as a hypersurface in the quasi-Euclidean space (6.1)
\be
\cS_r\subset\cQ_N,
\ee
we can arrive at

\ses\ses

{\large Theorem} 2.17. {\it The $\cQ_N$-induced geometry on the sphere}
(6.32) {\it is  a geometry of the constant curvature $h^2/r^2$, where
$h$ is the parameter} (2.13).

\ses\ses

This Theorem states that for the curvature $\cC$ of the hypersurface
(6.32)
we should get
\be
\cC_{quasi-Euclidean}=
h^2\cC_{Euclidean},
\ee
where
\be
\cC_{Euclidean}=\fr1{r^2}.
\ee
To verify this,
it is sufficient to note that for any admissible parameterizations
\be
t^p=t^p(u^a)
\ee
of $\cS_r$ the projection factors
\be
T_a^p=\D{t^p}{u^a},
\ee
where $a=1,\dots,N-1$,
satisfy the identity
\be
L_pT_a^p\equiv0
\ee
and, therefore,
the $L_p$-factors  disappear
in the induced metric tensor
\be
q_{ab}=
T_a^pT_b^qn_{pq}
\equiv
T_a^pT_b^q\Bigl(
\fr1{h^2}r_{pq}-\fr14G^2L_pL_q
\Bigr)
\ee
(the formula (5.49) has been applied) leaving us with the represention
\be
q_{ab}=T_a^pT_b^qr_{pq}/h^2.
\ee


A particular convenient parameterizations way is
$$
t^a=u^a, \qquad t^N=\sqrt{r^2-u^2}
$$
with
$$
u=\sqrt{r_{ab}u^au^b},
$$
entailing
$$
T_a^N=-\fr{u_a}{\sqrt{r^2-u^2}}, \qquad T_b^a=\de_b^a
$$
together with
$$
h^2q_{ab}(u)=r_{ab}+\fr{u_au_b}{r^2-u^2},
$$
\ses
$$
\fr1{h^2}q^{ab}(u)=r^{ab}-\fr 1{r^2}u^au^b,
$$
and
$$
\det(q_{ab}(u))=h^{2(1-N)}
\fr{r^2}{r^2-u^2}\det(r_{ab}).
$$
\ses\\
Constructing the associated Christoffel symbols yields
$$
I_{ab,e}=
\fr{r_{ab}u_e}{r^2-u^2}\fr1{h^2}
+
\fr{u_au_bu_e}{(r^2-u^2)^2}\fr1{h^2}
=\fr{u_e}{r^2-u^2}q_{ab},
\qquad I_a{}^e{}_b=\fr{h^2}{r^2}u^eq_{ab},
$$
which entails
that the associated  curvature tensor
$$
R_e{}^c{}_{ab}(u) :=\D{I_e{}^c{}_a}{u^b}-\D{I_e{}^c{}_b}{u^a}+
I_e{}^d{}_aI_d{}^c{}_b-I_e{}^d{}_bI_d{}^c{}_a
$$
of the sphere (6.33) proves entailing the simple form
$$
R_{ecab}(u)=h^2\fr1{r^2}(q_{ea}q_{cb}-q_{eb}q_{ac})
$$
that directly manifests  validity of the equality (6.34).

\bc
{\bf Chapter 3: Quasi-Euclidean consideration}
\ec

\ses
\ses

For the space under study, the geodesics should be obtained as
solutions to the equations (2.88)-(2.89) of Chapter 2 through
well-known arguments. To avoid complications of calculations
involved, it proves convenient to transfer preliminary the
consideration into the quasi-Euclidean approach (exposed in
Section 2.5). Below, we first consider in Section 3.1 the geodesic
equation exploiting attentively the fact that the respective
Christoffel symbols
 $\3Npqr$
are of sufficiently simple structure.
Surprisingly, the equation admits a simple and
explicit general solution, as this will be shown in great detail.
 After that, the angle between two vectors is explicated.
The remarkable result is that
the
angle $\al$ found
 is a factor of the Euclidean one,
the angle being normalized such that the Cosine Theorem of
ordinary form be rigorously valid if the Euclidean angle is
replaced by the $\al$. The respective scalar product ensues. In
Section 3.2, we introduce the associated two-vector metric tensor
and demonstrate that at  equality of vectors the tensor reduces
exactly to the one-vector FMT of the $\cE_g^{PD}$-space. The
orthonormal frame thereto is also obtained in a lucid explicit
form. After that, in Section 3.3, the possibility of converting
the theory into the co-approach is presented, and in Section 3.4
the $\cE_g^{PD}$-extension of the parallelogram law of vector
addition is derived; it occurs also possible to find the
respective sum vector and the difference vector in a nearest
approximation with respect to the characteristic parameter $g$.

 \ses \ses

\setcounter{sctn}{1}
\setcounter{equation}{0}
3.1. {\it
Derivation of geodesics
and angle
in
associated
 quasi-Euclidean space.
}
We start with searching for a general solution to the quasi-Euclidean
geodesic equation
which, in terms of
the coefficients
$\3Npqr$
(given by  Eq. (6.16) in Chapter 2),
reads
\be
\fr
{d^2t^p}
{ds^2} +\3Nqpr(g;t)
\fr{dt^q}{ds}
\fr{dt^r}{ds}
=0.
\ee
Accordingly, we put
\be
\sqrt{
g_{pq}(g;R)dR^pdR^q}
=
\sqrt{
n_{pq}(g;t)dt^pdt^q}
\ee
and
\be
R^p(s)=
\mu^p(g;t^r(s))
\ee
together with
\be
\fr
{dR^p(s)}
{ds}
=
\mu^p_q(g;t^r(s))
\fr
{dt^q(s)}
{ds},
\ee
where
$
\mu^p(g;t^r)
$
and
$
\mu^p_q(g;t^r)
$
are the coefficients that were given  in Chapter 2
by Eqs. (5.14)-(5.15) and
(5.38)-(5.40), respectively.
Let a curve $C$:\,
$t^p=t^p(s)$
be considered,
with
the
{\it
arc-length parameter} $s$ along the curve
being defined by the help of the differential
\be
ds=\sqrt{
n_{pq}(g;t)dt^pdt^q},
\ee
where
$
n_{pq}(g;t)
$
is the associated quasi-Euclidean metric tensor
(Eq. (5.49) in  Chapter 2).
Respectively,
the
{\it
 tangent vectors
}
\be
u^p=\fr{dt^p}{ds}
\ee
to the curve
are unit,
 in the sense that
\be
n_{pq}(g;t)u^pu^q=1.
\ee
Since
$L_p=\partial S/\partial t^p$,
we have
\be
L_pu^p=\fr{dS}{ds}.
\ee
Here,
$
S^2(t)=
n_{pq}(g;t)t^pt^q
=
r_{pq}t^pt^q
$
(see Eqs. (5.3) and (5.46) in  Chapter 2).
Using (1.1)
leads
to the following
equation for geodesics in the quasi-Euclidean space:
\be
\fr{d^2{\bf t}}{ds^2}=
\fr14G^2\fr{{\bf t}}{S^2}H_{pq}u^pu^q,
\ee
\ses
where
$H_{pq}=h^2(n_{pq}-L_pL_q)
$
(see Eq. (6.4) in  Chapter 2)
and
$
{\bf t}=\{t^p\}.
$
We obtain
\be
\fr{d^2{\bf t}}{ds^2}=
\fr14g^2\fr{{\bf t}}{S^2}
\lf(
1-(\fr
{dS}{ds})^2
\rg)
=
\fr14g^2(a^2-b^2)\fr{{\bf t}}{S^4}
\ee
\ses
and
\be
\fr{d^2{\bf t}}{ds^2}=
\fr14g^2(a^2-b^2)\fr{{\bf t}}{S^4}
\ee
with
\be
S^2(s)=a^2+2bs+s^2,
\ee
where $a$ and $b$ are two constants of integration.
The formula (1.12) is valid because from
\be
\fr12\fr{ dS^2}{ ds}=
r_{pq}t^pu^q
\ee
one can deduce
\ses\\
$$
\fr12\fr{ dS^2}{ ds^2}=
r_{pq}u^pu^q
+t_q\fr{ du^q}{ ds}
=(h^2
n_{pq}+\fr14g^2L_pL_q)u^pu^q+
\fr14g^2\lf(1-(\fr{ dS}{ ds})^2\rg)
$$
\ses
\ses
$$
=h^2+\fr14g^2
(\fr{ dS}{ ds})^2+
\fr14g^2\lf(1-(\fr{ dS}{ ds})^2\rg)
=1.
$$

If we put \be S (\Delta s)= \sqrt{ a^2+2b \Delta s + (\Delta s)^2
} \ee and \be {\bf t}_1 = {\bf t}(0), \qquad {\bf t}_2 = {\bf t}
(\Delta s), \ee then we get \be a= \sqrt{ ( {\bf t}_1 {\bf t}_1 )
} \ee and \be S (\Delta s)= \sqrt{ ({\bf t}_2 {\bf t}_2 ) }. \ee
Here, $ {\bf t}_1 $ and $ {\bf t}_2 $ are two vectors emanated
from  the fixed origin ``O"; they point to the beginning of the
geodesic and to the end of the geodesic, respectively. The
parenthesis couple $ (..) $ is used to denote the Euclidean scalar
product, so that $ ( {\bf t}_1 {\bf t}_1 ) = r_{pq}t_1^pt_1^q$, $
( {\bf t}_2 {\bf t}_2 ) = r_{pq}t_2^pt_2^q$, $ ( {\bf t}_1 {\bf
t}_2 ) = r_{pq}t_1^pt_2^q$; $r_{pq}=\de_{pq}$ in case of
orthogonal basis ($\de$ stands for the Kronecker symbol).

What value should be prescribed to the scalar product
$
(
{\bf t}_1
{\bf t}_2
)
$? A special analysis shows that the respective correct choice compatible
with the
geodesic equation solution (see Eqs. (1.41) and (1.42) below) should read
\be
(
{\bf t}_1
{\bf t}_2
)
=
 a
S
(\Delta s)
\cos
\Bigl[
h
\arctan\fr{
\sqrt{a^2-b^2}
\,
\Delta s
}{a^2+b
\Delta s
}
\Bigr].
\ee

From (1.16)-(1.18) it directly follows that
\be
\fr{
\sqrt{a^2-b^2}
\,
\Delta s
}
{a^2+b
\Delta s
}
=
\tan
\Bigl[
\fr1h\arccos
\fr{
(
{\bf t}_1
{\bf t}_2
)}
{
\sqrt{(
{\bf t}_1
{\bf t}_1
)}
\,
\sqrt{(
{\bf t}_2
{\bf t}_2
)}
}
\Bigr],
\ee
which entails
$$
\cos
\Bigl[
h\arctan\fr{
\sqrt{a^2-b^2}
\,
\Delta s
}{a^2+b
\Delta s
}
\Bigr]
=
\fr{
(
{\bf t}_1
{\bf t}_2
)}
{
\sqrt{(
{\bf t}_1
{\bf t}_1
)}
\,
\sqrt{(
{\bf t}_2
{\bf t}_2
)}
}
$$
and
$$
\sin \Bigl[ h\arctan\fr{ \sqrt{a^2-b^2} \, \Delta s }{a^2+b \Delta
s } \Bigr] = \fr { u ( {\bf t}_1, {\bf t}_2 ) }
 { \sqrt{( {\bf t}_1 {\bf t}_1 )} \,
\sqrt{( {\bf t}_2 {\bf t}_2 )} },
$$
where \be
 u ( {\bf t}_1, {\bf t}_2 )=\sqrt{ ( {\bf t}_1 {\bf t}_1
) ( {\bf t}_2 {\bf t}_2 ) - ( {\bf t}_1 {\bf t}_2 )^2 }.
 \ee
 The above  equalities   suggest the
idea to introduce

\ses
\ses

{\large Definition}.
{\it
The
$\cE_g^{PD}$-associated angle
}
is given by
\be
\al
(
{\bf t}_1,
{\bf t}_2
)
:=
\fr1h\arccos
\fr{
(
{\bf t}_1
{\bf t}_2
)}
{
\sqrt{(
{\bf t}_1
{\bf t}_1
)}
\,
\sqrt{(
{\bf t}_2
{\bf t}_2
)}
},
\ee
 so that
\be
\al
=\fr1h\al_{Euclidean}.
\ee

\ses
\ses

Such an angle is obviously
{\it
additive}:
\be
\al(
{\bf t}_1,
{\bf t}_3
)
=
\al(
{\bf t}_1,
{\bf t}_2
)
+
\al
(
{\bf t}_2,
{\bf t}_3
)
\ee
in planar case of the vector triple $\{
{\bf t}_1,
{\bf t}_2,
{\bf t}_3
\}$.

The traditional vanishing at equal vectors holds:
\be
\al(
{\bf t},
{\bf t}
)=0.
\ee

With any admissible prescribed parameter value $g$,
the range
\be
0\le
\al(
{\bf t}_1,
{\bf t}_2
)\le\fr1h\pi
\ee
corresponds to
 the canonical Euclidean range $(0,\pi)$.

With the angle (1.21), we ought to propose

\ses
\ses

{\large Definition}.
Given two vectors
$
{\bf t}
$
and
$
{\bf t}^{\perp},
$
 the vectors are said to be
{\it
$\cE_g^{PD}$-perpendicular,
}
if
\be
\cos
\lf(
\al
(
{\bf t},
{\bf t}
^{\perp})
\rg)
=0.
\ee

Since the vanishing (1.25) implies \be \al ( {\bf t}, {\bf t}
^{\perp}) = \fr{\pi}2, \ee in view of (1.22) we ought to conclude
that \be \al_{Euclidean} ( {\bf t}, {\bf t} ^{\perp}) =
\fr{\pi}2h\le \fr{\pi}2. \ee Therefore, vectors perpendicular in
the quasi-Euclidean sense proper look like acute vectors as
observed from associated Euclidean standpoint proper.

With the equality
\be
(
\sqrt{a^2-b^2}
\,
\Delta s)^2
+(a^2+b\Delta s)^2\equiv
a^2S^2(\Delta s)
\ee
(an implication of Eq. (1.14))
we
also establish
the relations
\be
\sqrt{a^2-b^2}
\,
\Delta s=aS(\Delta s)
\sin
\al
\ee
and
\ses\\
\be
a^2+b\Delta s=
aS(\Delta s)
\cos
\al,
\ee
where $\al$ is exactly the angle (1.21).
They entail the equalities
\be
\fr b{\sqrt{a^2-b^2}}
=
\fr{
S(\Delta s)
\cos
\al
-a}{
S(\Delta s)
\sin
\al
},
\qquad
\fr{b^2}{a^2}=1-\lf(
\fr{S(\De s)}{\De s}\rg)^2\sin^2\al
\ee
from which the quantity $b$ also can be explicated.

Thus
{\it
each member of the involved set
$\{
a, b, \Delta s, S(\Delta s)\}
$
can be explicitly expressed through the input vectors
$
{\bf t}_1
$
and
$
{\bf t}_2.
$
}
For various particular cases it is worth rewriting the equality (1.29) as
\be
S^2(\Delta s)=
(\Delta s)^2-a^2+
2(a^2+b
\Delta s
).
\ee

Thus we may naturally set forth  the following fundamental items:
\ses
\ses\\
\nin
{\it
The
quasi-Euclidean
Cosine Theorem
}
\ses\\
\be
(\Delta s)^2
=
S^2(\Delta s)
+a^2-
2aS(\Delta s)
\cos
\al.
\ee

\ses
\ses

\nin
{\it
The
quasi-Euclidean
Two-Point Length
}
\ses\\
\be
(\Delta s)^2
=
({\bf t}_1
{\bf t}_1
)
+
({\bf t}_2
{\bf t}_2
)
-2
\sqrt
{
({\bf t}_1
{\bf t}_1
)
}
\,
\sqrt
{
({\bf t}_2
{\bf t}_2
)
}
\cos
\al.
\ee

\ses
\ses

\nin
{\it
The
quasi-Euclidean
Scalar Product
}
\ses\\
\be
<{\bf t}_1,
{\bf t}_2
>
=
\sqrt
{
({\bf t}_1
{\bf t}_1
)
}
\,
\sqrt
{
({\bf t}_2
{\bf t}_2
)
}
\cos
\al.
\ee

\ses
\ses

\nin
{\it
The
quasi-Euclidean
Perpendicularity
}
\ses\\
\be
<{\bf t},
{\bf t}^{\perp}
>
=
\sqrt
{
({\bf t}{\bf t})
}
\,
\sqrt
{
({\bf t}^{\perp}
{\bf t}^{\perp})
}.
\ee

\ses
\ses

The identification
\ses\\
\be
|
{\bf t}_2
\ominus
{\bf t}_1
|^2
=(\Delta s)^2
\ee
yields another lucid representation
\be
|
{\bf t}_2
\ominus
{\bf t}_1
|^2
=
({\bf t}_1
{\bf t}_1
)
+
({\bf t}_2
{\bf t}_2
)
-2
\sqrt
{
({\bf t}_1
{\bf t}_1
)
}
\,
\sqrt
{
({\bf t}_2
{\bf t}_2
)
}
\cos
\al\,.
\ee
Here the symmetry holds:
\be
|
{\bf t}_2
\ominus
{\bf t}_1
|
=
|
{\bf t}_1
\ominus
{\bf t}_2
|.
\ee

Comparing
$$
({\bf t}_1{\bf t}_2)=
\sqrt{({\bf t}_1{\bf t}_1)}\,\sqrt{({\bf t}_2{\bf t}_2)}
\cos
h\al
=
\sqrt{({\bf t}_1{\bf t}_1)}\,\sqrt{({\bf t}_2{\bf t}_2)}
\cos
\al_{Euclidean}
$$
(see Eqs. (1.18)-(1.21)) with (1.36) shows that
$$
<{\bf t}_1,{\bf t_2}>\ne
({\bf t}_1{\bf t_2}) \quad {\rm unless}~ g=0.
$$

The consideration can be completed by

\ses
\ses

{\large Proposition} 3.1.
{\it
A general solution to the geodesic equation {\rm (1.11)}
can explicitly be found as follows:
\ses
$$
{\bf t}(s)=
$$
\ses
\ses
\be
=
\fr {S(s)}{
 a
}
\fr
{\sin
\Bigl[
h\arctan\fr{
\sqrt{a^2-b^2}
\,
(\Delta s-s)
}{a^2+b
\Delta s
+(b+
\Delta s
)
 s
}
\Bigr]
}{
\sin
\Bigl[
h\arctan\fr{
\sqrt{a^2-b^2}
\,
\Delta s
}{a^2+b
\Delta s
}
\Bigr]
}
\,
{\bf t}_1
+
\fr{ S(s)}{
S
(\Delta s)
}
\fr
{\sin
\Bigl[
h\arctan\fr{
\sqrt{a^2-b^2}
\,
s
}{a^2+bs
}
\Bigr]
}{
\sin
\Bigl[
h\arctan\fr{
\sqrt{a^2-b^2}
\,
\Delta s
}{a^2+b
\Delta s
}
\Bigr]
}
\,
{\bf t}_2.
\ee
}

\ses
\ses

The validity of this proposition can be processed by direct insertion of
the expression (1.41)
in (1.11).

\ses

From (1.41) the equality
\ses
\be
(
{\bf t}(s)
{\bf t}(s)
)
=S^2(s)
\ee
follows, in agreement with (1.12).
Also, it is useful to note that
$$
\arctan\fr{
\sqrt{a^2-b^2}
\,
(\Delta s-s)
}{a^2+b
\Delta s
+(b+
\Delta s
)
 s
}
+
\arctan\fr{
\sqrt{a^2-b^2}
\,
s
}{a^2+bs
}
=
\arctan\fr{
\sqrt{a^2-b^2}
\,
\Delta s
}{a^2+b
\Delta s
},
$$
where the right-hand part does not involve the variable $s$. Using
(1.19) and (1.21), we conclude \be \arctan\fr{ \sqrt{a^2-b^2} \,
(\Delta s-s) }{a^2+b \Delta s +(b+ \Delta s )
 s
}  =\al-\nu, \ee where \be
 \nu=\arctan\fr{ \sqrt{a^2-b^2} \, s
}{a^2+bs }. \ee

In terms of this angle, the solution (1.41) reads simply
 \be {\bf
t}(s)= \fr{S(s)}a \fr{\sin(h(\al-\nu))}{\sin(h\al)} {\bf t}_1 +
\fr{S(s)}{S(\De s)} \fr{\sin(h\nu)}{\sin(h\al)} {\bf t}_2.
 \ee

 The
identification
$$
{{\bf t}(s)}
_{\Bigl|_
{s=0}
\Bigr{.}}
=
{\bf t}_1,
\qquad
{{\bf t}(s)}
_{\Bigl|_
{s=\Delta s}
\Bigr{.}}
=
{\bf t}_2.
$$
can readily be verified.
The Euclidean limit proper for the solution (1.41) is
\ses\\
$$
{{\bf t}(s)}
_{\Bigl|_
{g=0}
\Bigr{.}}
=
\fr{
(\Delta s-s)
{\bf t}_1
+
s{\bf t}_2
}{
\Delta s
}
=
{\bf t}_1+
(
{\bf t}_2-
{\bf t}_1
)\fr s{\Delta s},
$$
so that the geodesics are got simplified to be straight lines.

Since the general solution (1.41) is such that the right-hand side
is spanned by two fixed vectors,
$
{\bf t}_1
$
and
$
{\bf t}_2,
$
we are entitled concluding that

\ses
\ses

{\large Proposition} 3.2.
 {\it
Against the quasi-Euclidean treatment,
the geodesics under study are plane curves.
}

\ses
\ses

Calculating the first derivative
$$
{\bf v}(s) := \fr{d{\bf t} (s)} {ds}
$$
of (1.41) yields
the formula
\ses\\
\be
{\bf v}(s) = \fr{b+s}{S^2(s)} {\bf t}(s) \, - \fr{
\sqrt{a^2-b^2} \, h } { aS(s)} \fr {\cos (h(\al-\nu)) }{ \sin (h\al)} \, {\bf t}_1
+
 \fr{ \sqrt{a^2-b^2} \, h}{ S(s) S (\Delta s) } \,
 \fr {\cos (h\nu) }{ \sin (h\al)} \, {\bf t}_2.
 \ee
  The right-hand part here is such that
$$
{\bf t}(s)
\lf(
{\bf v}(s)
-
\fr{b+s}{S^2(s)}
{\bf t}
(s)\rg)
=0.
$$
From the latter  observation we get the  equality
$$
{\bf t}(s)
{\bf v}
(s)
=b+s
$$
 which is tantamount to (1.12).

Also, the {\it initial velocity}
\ses
$$
{\bf v}_1 :=
\fr{d{\bf t}
}
{ds}
_{\Bigl|_
{
s
=
0}
\Bigr{.}}
$$
and the {\it final velocity}
 $$ {\bf v}_2 := \fr{d{\bf t} }
{ds} _{\Bigl|_ { s = \Delta s } \Bigr{.}} $$ are found from (1.46)
to be \ses
\ses\\
$$
{\bf v}_1 = \fr{b}{a^2} {\bf t}_1 - \fr { \sqrt{a^2-b^2} \, h }
{a^2} \fr { \cos(h\al) } { \sin (h\al) } \, {\bf t}_1 + \fr {
\sqrt{a^2-b^2} \, h } {aS(\Delta s)} \fr1 { \sin(h\al) } \, {\bf
t}_2
$$
\ses
\ses
\ses
\ses
\ses
\ses
$$
=
\fr{b}{a^2}
{\bf t}_1
-
\fr
{
\sqrt{a^2-b^2}
\, h
}
{a^2}
\fr
{
(
{\bf t}_1
{\bf t}_2
)
}
{
u({\bf t}_1,{\bf t}_2)
}
\,
{\bf t}_1
+
\sqrt{a^2-b^2}
\, h
\fr
{
1
}
{
u({\bf t}_1,{\bf t}_2)
}
\,
{\bf t}_2
$$
\ses\\
and
\ses\\
$$
{\bf v}_2 = \fr{b+\Delta s}{S^2(\Delta s)} {\bf t}_2 - \fr {
\sqrt{a^2-b^2} \, h } {aS(\Delta s)} \fr1 { \sin (h\al) } \, {\bf
t}_1 + \fr { \sqrt{a^2-b^2} \, h } {S^2(\Delta s)} \fr { \cos
(h\al) } { \sin (h\al) } \, {\bf t}_2
$$
\ses
\ses
\ses
\ses
\ses
\ses
\ses
$$
=
\fr{b+\Delta s}{S^2(\Delta s)}
{\bf t}_2
-
\sqrt{a^2-b^2}
\, h
\fr
{
1
}
{
u({\bf t}_1,{\bf t}_2)
}
\,
{\bf t}_1
+
\fr
{
\sqrt{a^2-b^2}
\, h
}
{S^2(\Delta s)}
\fr
{
(
{\bf t}_1
{\bf t}_2
)
}
{
u({\bf t}_1,{\bf t}_2)
}
\,
{\bf t}_2
$$
(with having used (1.19)), where $ u({\bf t}_1,{\bf t}_2) $ is the
function (1.20).

Notice also that
$$
(
{\bf t}_1
{\bf v}_1
)=b,
\qquad
(
{\bf t}_2
{\bf v}_2
)=b+\De s,
$$
\ses
$$
( {\bf v}_1 {\bf v}_2 )= \lf(1+\fr{G^2}4\fr{b^2}{a^2}\rg)h^2
\equiv\fr{(\De s)^2h^2+(1-h^2) ( {\bf t}_2 {\bf
t}_2)\sin^2\al}{(\De s)^2}
$$ and
$$
n_{pq}(g;
{\bf t}_1)
v_1^pv_1^q=1.
$$

\ses

The difference between the vectors $ {\bf v}_2 $ and $ {\bf v}_1 $
 can be found to read
\ses
$$
{\bf v}_2
-
{\bf v}_1
=
\Bigl[
\fr{b+\Delta s}{S^2(\Delta s)}
+
\fr
{
\sqrt{a^2-b^2}
\, h
}
{S^2(\Delta s)}
\fr
{
(
{\bf t}_1
{\bf t}_2
)
}
{
u({\bf t}_1,{\bf t}_2)
}
-
\sqrt{a^2-b^2}
\, h
\fr
{
1
}
{
u({\bf t}_1,{\bf t}_2)
}
\Bigr]
{\bf t}_2
$$
\ses
\ses
\ses
$$
-
\Bigl[
\fr{b}{a^2}
+
\sqrt{a^2-b^2}
\, h
\fr1
{
u({\bf t}_1,{\bf t}_2)
}
-
\fr
{
\sqrt{a^2-b^2}
\, h
}
{a^2}
\fr
{
(
{\bf t}_1
{\bf t}_2
)
}
{
u({\bf t}_1,{\bf t}_2)
}
\Bigr]
{\bf t}_1.
$$
\ses
It is useful to verify that
$$
{{\bf v}_2}_{\Bigl|_ { g = 0} \Bigr{.}} = {{\bf v}_1}_{\Bigl|_ { g
= 0} \Bigr{.}} = \fr{ {\bf t}_2- {\bf t}_1 } {\Delta s}
$$
in compliance with the Euclidean rule proper.

\ses

Other convenient representations for these vectors are
$$
 {\bf v}_1\De s=
(-1+\fr1{({\bf t}_1{\bf t}_1)}C_1){\bf t}_1+
h \fr{\sin\al}{\sin\al_{Euclidean}}{\bf t}_2,
 $$
 \ses
 $$
  {\bf v}_2\De s= (1-\fr1{({\bf t}_2{\bf t}_2)}C_1){\bf t}_2-
h \fr{\sin\al}{\sin\al_{Euclidean}}{\bf t}_1,
$$
  and
$$
({\bf v}_2 -{\bf v}_1)\De s=\lf(1-h\fr{\sin\al}{\sin\al_{Euclidean}}
\rg)({\bf t}_2
+
{\bf t}_1
)-C_1(
\fr{
{\bf t}_2}{({\bf t}_2{\bf t}_2)}+
\fr{
{\bf t}_1}{({\bf t}_1{\bf t}_1)}
),
$$
where
$$
C_1=
\sqrt{({\bf t}_1{\bf t}_1)}\sqrt{({\bf t}_2{\bf t}_2)}
\cos\al-h
({\bf t}_1{\bf t}_2)
\fr{\sin\al}{\sin\al_{Euclidean}}.
$$

\ses \ses

Using above formulae, we can turn the solution (1.40) into the
{\it initial-date form}:

\ses \ses

{\large Proposition} 3.3. {\it The general initial-date solution
to the geodesic equation} (1.1) {\it reads \be {\bf t}(s)=m(s){\bf
t}_1+n(s){\bf v}_1 \ee with
$$
m(s)= - \fr {b S(s)} { a\sqrt{a^2-b^2} \, h} \sin(h\al) + \fr{
S(s) } { a} \cos (h\al) + \fr {S(s)}{
 a
} \fr {\sin(h(\al-\nu))
 }{
\sin(h\al)  }
$$
\ses and
$$
n(s)= \fr {a S(s)} { \sqrt{a^2-b^2} \, h} \sin (h\nu).
$$
}

\ses \ses

Below the picture symbolizes the character of the solution (1.47).

\ses

\begin{figure}[!ht]
\vspace{2cm} \centering
\includegraphics[width=14cm]{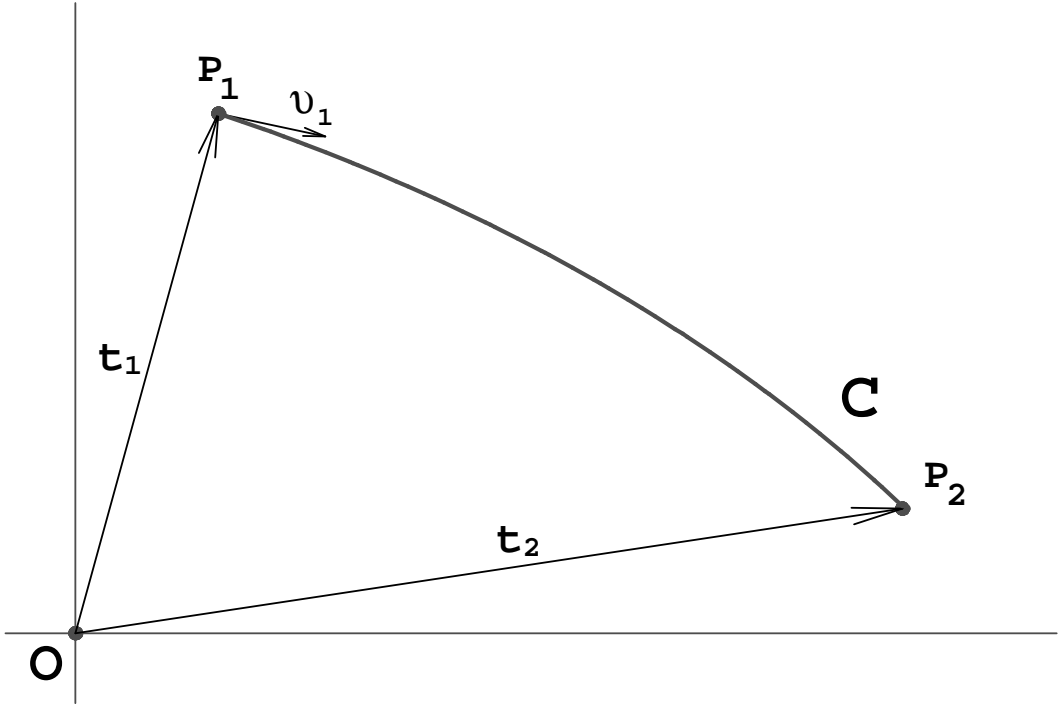}
\caption{[The geodesic $C$ and the initial velocity ${\bf v}_1$]}
\end{figure}

\clearpage

For the {\it derivative vectors}
$$
{\bf b}_1
:=
\fr12
\D{
|{\bf t}_2\ominus
{\bf t}_1|^2
}
{
{\bf t}_1
},
\qquad
{\bf b}_2
:=
\fr12
\D{
|{\bf t}_2\ominus
{\bf t}_1|^2
}
{
{\bf t}_2
},
$$
we can obtain the simple representations
\ses
\ses\\
$$
{\bf b}_1
=
{\bf t}_1
-
\fr
{
{\bf t}_1
}
{
\sqrt
{
(
{\bf t}_1
{\bf t}_1
)
}
}
\sqrt
{
(
{\bf t}_2
{\bf t}_2
)
}\,
\cos
\al
-
\fr{
\sqrt{(
{\bf t}_2
{\bf t}_2
)}
}
{
h
\sqrt{(
{\bf t}_1
{\bf t}_1
)}
}
{\bf d}_1
\sin
\al
$$
and
$$
{\bf b}_2
=
{\bf t}_2
-
\fr
{
{\bf t}_2
}
{
\sqrt
{
(
{\bf t}_2
{\bf t}_2
)
}
}
\sqrt
{
(
{\bf t}_1
{\bf t}_1
)
}\,
\cos
\al
-
\fr{
\sqrt{(
{\bf t}_1
{\bf t}_1
)}
}
{
h
\sqrt{(
{\bf t}_2
{\bf t}_2
)}
}
{\bf d}_2
\sin
\al,
$$
\ses
where
the convenient vectors
\ses\\
$$
{\bf d}_1
=
\fr{
(
{\bf t}_1
{\bf t}_1
)
{\bf t}_2
-
(
{\bf t}_1
{\bf t}_2
)
{\bf t}_1
}
{
u({\bf t}_1,{\bf t}_2)
},
\qquad
{\bf d}_2
=
\fr{
(
{\bf t}_2
{\bf t}_2
)
{\bf t}_1
-
(
{\bf t}_1
{\bf t}_2
)
{\bf t}_2
}
{
u({\bf t}_1,{\bf t}_2)
}
$$
have been introduced.
It can readily be verified that
$$
(
{\bf t}_1
{\bf d}_1
)=0,
\qquad
(
{\bf t}_2
{\bf d}_2
)=0,
$$
\ses
\ses
$$
(
{\bf d}_1
{\bf d}_2
)=-
(
{\bf t}_1
{\bf t}_2
),
\qquad
(
{\bf d}_1
{\bf d}_1
)=
(
{\bf t}_1
{\bf t}_1
),
\qquad
(
{\bf d}_2
{\bf d}_2
)=
(
{\bf t}_2
{\bf t}_2
),
$$
\ses
\ses
$$
(
{\bf d}_1
{\bf t}_2
)
=
(
{\bf t}_1
{\bf d}_2
)
=
u({\bf t}_1,{\bf t}_2),
$$
\ses
and
\ses\\
$$
{\bf t}_1
{\bf b}_1
+
{\bf t}_2
{\bf b}_2
=
2
|{\bf t}_2\ominus
{\bf t}_1|^2,
$$
\ses
together with
$$
\lim_{
{\bf t}_2\to
{\bf t}_1
}
\Bigl\{
{\bf b}_1
\Bigl\}=
\lim_{
{\bf t}_2\to
{\bf t}_1
}
\Bigl\{
{\bf b}_2
\Bigl\}=0.
$$

For the respective products of  vectors  we obtain
\ses
$$
(
{\bf b}_1
{\bf b}_1
)
=
(
{\bf t}_1
{\bf t}_1
)
+
(
{\bf t}_2
{\bf t}_2
)
-2
\sqrt
{
(
{\bf t}_1
{\bf t}_1
)
}
\sqrt
{
(
{\bf t}_2
{\bf t}_2
)
}
\,
\cos
\al
+
\lf(\fr1{h^2}-1\rg)
(
{\bf t}_2
{\bf t}_2
)
\sin^2
\al,
$$
\ses
\ses
\ses
$$
(
{\bf b}_2
{\bf b}_2
)
=
(
{\bf t}_2
{\bf t}_2
)
+
(
{\bf t}_1
{\bf t}_1
)
-2
\sqrt
{
(
{\bf t}_2
{\bf t}_2
)
}
\sqrt
{
(
{\bf t}_1
{\bf t}_1
)
}
\,
\cos
\al
+
\lf(\fr1{h^2}-1\rg)
(
{\bf t}_1
{\bf t}_1
)
\sin^2
\al,
$$
\ses
and
\ses
$$
(
{\bf b}_1
{\bf b}_2
)
=
-
\Bigl[
\lf(
\fr
{
\sqrt{
(
{\bf t}_1
{\bf t}_1
)
}
}
{
\sqrt
{
(
{\bf t}_2
{\bf t}_2
)
}
}
+
\fr
{
\sqrt{
(
{\bf t}_2
{\bf t}_2
)
}
}
{
\sqrt
{
(
{\bf t}_1
{\bf t}_1
)
}
}
-2\cos\al
\rg)
\cos
\al
+
\lf(\fr1{h^2}-1\rg)\sin^2\al
\Bigr]
(
{\bf t}_1
{\bf t}_2
)
$$
\ses
\ses
$$
-
\fr1h
\lf(
\fr
{
\sqrt{
(
{\bf t}_1
{\bf t}_1
)
}
}
{
\sqrt
{
(
{\bf t}_2
{\bf t}_2
)
}
}
+
\fr
{
\sqrt{
(
{\bf t}_2
{\bf t}_2
)
}
}
{
\sqrt
{
(
{\bf t}_1
{\bf t}_1
)
}
}
-2\cos\al
\rg)
u({\bf t}_1,{\bf t}_2)
\sin
\al.
$$
\ses\\
The vectors ${\bf b}_1$ and ${\bf b}_2$ extend the difference vectors,
for in the  Euclidean limit we have simply
$$
{{\bf b}_1}
_{\Bigl|_
{
g
=
0}
\Bigr{.}}
=
{\bf t}_1
-
{\bf t}_2,
\qquad
{{\bf b}_2}
_{\Bigl|_
{
g
=
0}
\Bigr{.}}
=
{\bf t}_2
-
{\bf t}_1.
$$

The following limit
\be
\lim_{
{\bf t}_2\to
{\bf t}_1
}
\Bigl\{
\fr{
(
{\bf t}_1
{\bf t}_1
)
(
{\bf t}_2
{\bf t}_2
)
}
{
h
\sqrt{(
{\bf t}_1
{\bf t}_1
)}
\sqrt{(
{\bf t}_2
{\bf t}_2
)}
}
\fr{
\sin
\Bigl[
\fr1h\arccos
\fr{
(
{\bf t}_1
{\bf t}_2
)}
{
\sqrt{(
{\bf t}_1
{\bf t}_1
)}
\,
\sqrt{(
{\bf t}_2
{\bf t}_2
)}
}
\Bigr]
}
{
u({\bf t}_1,{\bf t}_2)
}
\Bigl\}
=
\fr1{h^2}
\ee
is important to note.

\ses
\ses

\setcounter{sctn}{2}
\setcounter{equation}{0}
3.2. {\it
The two-vector metric tensor and frame
in quasi-Euclidean space}.
Now we are able to introduce
the
{\it
 quasi-Euclidean two-vector metric tensor}
$
n(g;
{\bf t}_1,
{\bf t}_2
)
$
by the components
\be
n_{pq}(g;
{\bf t}_1,
{\bf t}_2
) :=
\Dd{
<
{\bf t}_1,
{\bf t}_2>
}
{
t^q_2
}
{
 t^p_1
}
=
-\fr12
\Dd{
|{\bf t}_2\ominus
{\bf t}_1|^2
}
{
t^q_2
}
{
 t^p_1
}.
\ee
Straightforward calculations
(on the basis of (1.32) and (1.19))
show that
\ses\\
$$
n_{pq}(g;
{\bf t}_1,
{\bf t}_2
)
=
\fr{
(
{\bf t}_1
{\bf t}_1
)
(
{\bf t}_2
{\bf t}_2
)
}
{
h
\sqrt{(
{\bf t}_1
{\bf t}_1
)}
\sqrt{(
{\bf t}_2
{\bf t}_2
)}
}
\fr{
\sin
\al
}
{
u({\bf t}_1,{\bf t}_2)
}
r_{pq}
$$
\ses
\ses
\be
+
\fr{
1
}
{
\sqrt
{
(
{\bf t}_1
{\bf t}_1
)
}
\sqrt
{
(
{\bf t}_2
{\bf t}_2
)
}
}
A_1
 t_{1p}
t_{2q}
-
\fr{
1
}
{
h
\sqrt{(
{\bf t}_1
{\bf t}_1
)}
\sqrt{(
{\bf t}_2
{\bf t}_2
)}
}
A_2
d_{1p}
 d_{2q},
\ee
\ses\\
where
\be
A_1=
\cos
\al
-
\fr1h
(
{\bf t}_1
{\bf t}_2
)
\fr{
\sin
\al
}
{
u({\bf t}_1,{\bf t}_2)
}
\ee
and
\be
A_2=
\fr1h
\cos
\al
-
(
{\bf t}_1
{\bf t}_2
)
\fr{
\sin
\al
}
{
u({\bf t}_1,{\bf t}_2)
}.
\ee
\ses
\ses
For the determinant of the tensor (2.2) we find simply
\ses\\
\be
\det\lf(
n_{pq}(g;
{\bf t}_1,
{\bf t}_2
)
\rg)
=\lf(
\fr
{
\sqrt{({\bf t}_1{\bf t}_1)({\bf t}_2{\bf t}_2)}\sin\al
}
{
u({\bf t}_1,{\bf t}_2)
}
\rg)^{N-2}h^{-N}
\det\lf(r_{ab}\rg).
\ee

\ses
\ses

Owing to (1.48), we can establish
 the following fundamental identification:
\be
\lim_{
{\bf t}_2\to
{\bf t}_1=
{\bf t}
}
\Bigl\{
n_{pq}(g;
{\bf t}_1,
{\bf t}_2
)
\Bigl\}=
n_{pq}(g;t),
\ee
\ses
\ses
\ses
where
$
n_{pq}(g;
{\bf t}
)
$
is the quasi-Euclidean metric tensor (see (5.49) in Chapter 2).

Differentiating (2.2)  results in
$$
\D
{
n_{pq}(g;
{\bf t}_1,
{\bf t}_2
)
}
{
 t_1^s
}
=
-\fr1h
\fr
{
\sqrt
{
(
{\bf t}_2
{\bf t}_2
)
}
}
{
\sqrt
{
(
{\bf t}_1
{\bf t}_1
)
}
u
(
{\bf t}_1,
{\bf t}_2
)
}
A_2d_{1s}r_{pq}
$$
\ses
\ses
$$
+
\fr
1
{
\sqrt{
(
{\bf t}_1
{\bf t}_1
)
}
\sqrt{
(
{\bf t}_2
{\bf t}_2
)
}
}
A_1t_{2q}H_{sp}
(
{\bf t}_1
)
+
\fr1h
\fr
{
(
{\bf t}_1
{\bf t}_2
)
}
{
(
{\bf t}_1
{\bf t}_1
)
\sqrt{
(
{\bf t}_1
{\bf t}_1
)
}
\sqrt{
(
{\bf t}_2
{\bf t}_2
)
}
}
\fr1
{
u
(
{\bf t}_1,
{\bf t}_2
)
}
A_2d_{1s}t_{1p}t_{2q}
$$
\ses
\ses
$$
+
\fr1h
\fr
1
{
\sqrt{
(
{\bf t}_1
{\bf t}_1
)
}
\sqrt{
(
{\bf t}_2
{\bf t}_2
)
}
}
\fr1
{
u
(
{\bf t}_1,
{\bf t}_2
)
}
A_2
\times
$$
\ses
\ses
$$
\Bigl[
\fr1
{
(
{\bf t}_1
{\bf t}_1
)
}
\lf(
(
{\bf t}_1
{\bf t}_1
)
d_{2p}d_{2q}
+
(
{\bf t}_2
{\bf t}_2
)
d_{1p}d_{1q}
\rg)
d_{1s}
+
(
{\bf t}_1
{\bf t}_2
)
H_{ps}
(
{\bf t}_1
)
d_{2q}
-
(
{\bf t}_2
{\bf t}_2
)
H_{qs}
(
{\bf t}_1
)
d_{1p}
\Bigl]
$$
\ses
\ses
\be
+
\fr1h
\fr1
{
(
{\bf t}_1
{\bf t}_1
)
\sqrt{
(
{\bf t}_1
{\bf t}_1
)
}
\sqrt{
(
{\bf t}_2
{\bf t}_2
)
}
}
\Bigl[
(1-\fr1{h^2})\sin\al-
\fr
{
(
{\bf t}_1
{\bf t}_2
)
}
{u
(
{\bf t}_1,
{\bf t}_2
)
}
A_2
\Bigr]
d_{1s}d_{1p}d_{2q},
\ee
\ses\\
where  we have used the relations:
$$
 \D{\al} {t_1^s} = -\fr1h \fr
{ {\bf d}_1 } { ( {\bf t}_1 {\bf t}_1 ) }, \qquad \D{\fr1u}
{t_1^s} = -\fr1{u^2} {\bf d}_2,
$$ \ses \ses
$$ \D{A_1} {t_1^s}= \fr1h \fr { ( {\bf t}_1 {\bf t}_2 ) } { ( {\bf
t}_1 {\bf t}_1 ) u ( {\bf t}_1, {\bf t}_2 ) } A_2d_{1s}, $$
 \ses
\ses
$$
\D{A_2}
{t_1^s}=
\fr1h
\fr
{
(
{\bf t}_1
{\bf t}_2
)
}
{
(
{\bf t}_1
{\bf t}_1
)
u
(
{\bf t}_1,
{\bf t}_2
)
}
A_1
d_{1s}
-
(1-\fr1{h^2})
\fr
{
(
{\bf t}_2
{\bf t}_2
)
}
{
u
(
{\bf t}_1,
{\bf t}_2
)
}
\fr
{
\sin\al
}
{
u
(
{\bf t}_1,
{\bf t}_2
)
}
 d_{1s}
$$
\ses\ses
$$ = \fr{d_{1s}} { ( {\bf t}_1 {\bf t}_1 ) } \Bigl[
-(1-\fr1{h^2})\sin\al+ \fr { ( {\bf t}_1 {\bf t}_2 ) } {u ( {\bf
t}_1, {\bf t}_2 ) } A_2\Bigr].
$$ \ses Since $$ \lim_{ {\bf
t}_2\to {\bf t}_1 } \Bigl\{ {A_1} \Bigl\}= 1-\fr1{h^2}, \qquad
\lim_{ {\bf t}_2\to {\bf t}_1 } \Bigl\{ \fr {A_2} {u} \Bigl\}= 0,
$$ \ses \ses $$ \lim_{ {\bf t}_2\to {\bf t}_1= {\bf t} } \Bigl\{
\D { n_{pq}(g; {\bf t}_1, {\bf t}_2 ) } {
 t_1^s
} \Bigl\}= (1-\fr1{h^2}) \fr {t_q} { ( {\bf t} {\bf t} ) } H_{sp}(
{\bf t} ), $$ \ses and \ses \ses $$ \lim_{ {\bf t}_2\to {\bf t}_1=
{\bf t} } \Bigl\{ \D { n_{pq}(g; {\bf t}_1, {\bf t}_2 ) } {
 t_2^s
} \Bigl\}= (1-\fr1{h^2}) \fr {t_p} { ( {\bf t} {\bf t} ) } H_{sq}(
{\bf t} ), $$
 the fundamental consequence
\ses\\
$$ \lim_{ {\bf t}_2\to {\bf t}_1= {\bf t} } \Bigl\{ \D { n_{pq}(g;
{\bf t}_1, {\bf t}_2 ) } {
 t_1^s
}
+
\D
{
n_{pq}(g;
{\bf t}_1,
{\bf t}_2
)
}
{
 t_2^s
}
\Bigl\}=
\D
{
n_{pq}(g;
{\bf t}
)
}
{
 t^s
} $$
 is obtained.

The expansion with respect to an appropriate orthonormal frame $
f^R_p (g; {\bf t}_1, {\bf t}_2 ) $ can be found to read $$ n_{pq}
(g; {\bf t}_1, {\bf t}_2 ) = \sum_{R=1}^{N} f^R_p (g; {\bf t}_1,
{\bf t}_2 ) f^R_q (g; {\bf t}_2, {\bf t}_1 ) $$ \ses with
$$
\sqrt
{
h
\sqrt{(
{\bf t}_1
{\bf t}_1
)}
\sqrt{(
{\bf t}_2
{\bf t}_2
)}
}
f^R_p
(g;
{\bf t}_1,
{\bf t}_2
)
=
zh^R_p
$$
\ses
\ses
\ses
\ses
$$
-\fr1{
(
{\bf t}_1
{\bf t}_2
)
}
\Bigl[
z
-
\sqrt{z^2+
(
{\bf t}_1
{\bf t}_2
)
\Bigl(
h
\cos
\al
-
(
{\bf t}_1
{\bf t}_2
)
\fr{
\sin
\al
}
{
u({\bf t}_1,{\bf t}_2)
}
\Bigr)
}
\Bigr]
 t_2^R
 t_{1p}
$$
\ses \ses \ses \ses \ses $$ +\fr1{ ( {\bf t}_1 {\bf t}_2 )} \Bigl[
z-\sqrt{z^2+ ( {\bf t}_1 {\bf t}_2 ) \Bigl( \fr1h \cos \al - (
{\bf t}_1 {\bf t}_2 ) \fr{ \sin \al } { u({\bf t}_1,{\bf t}_2) }
\Bigr) } \Bigr]
 d_{2p}
d_1^R, $$ where
\ses\\
$$ z= \sqrt { ( {\bf t}_1 {\bf t}_1 ) ( {\bf t}_2 {\bf t}_2 ) \fr{
\sin \al } { u({\bf t}_1,{\bf t}_2) } }, $$ \ses \ses or finally,
\ses\\
$$
\sqrt
{
h
\sqrt{(
{\bf t}_1
{\bf t}_1
)}
\sqrt{(
{\bf t}_2
{\bf t}_2
)}
}
f^R_p
(g;
{\bf t}_1,
{\bf t}_2
)
$$
\ses
\ses
\ses
\ses
$$
=
zh^R_p
-\fr1{
(
{\bf t}_1
{\bf t}_2
)
}
\Bigl[
z-\sqrt{
h
(
{\bf t}_1
{\bf t}_2
)
\cos
\al
+
u({\bf t}_1,{\bf t}_2)
\sin
\al
}
\Bigr]
\,
 t_2^R
 t_{1p}
$$
\ses \ses \ses \ses \ses \ses \ses \ses $$ +\fr1{ ( {\bf t}_1 {\bf
t}_2 )} \Bigl[ z-\sqrt{ \fr1h ( {\bf t}_1 {\bf t}_2 ) \cos \al +
u({\bf t}_1,{\bf t}_2) \sin \al } \Bigr] \,
 d_{2p}
d_1^R. $$

Contracting the frame by vectors yields
\ses\\
$$
f^R_p
(g;
{\bf t}_1,
{\bf t}_2
)
t_1^p
$$
\ses
\ses
$$
=
\fr
1
{
\sqrt
{
h
\sqrt{
(
{\bf t}_1
{\bf t}_1
)}
\sqrt{
(
{\bf t}_1
{\bf t}_2
)}
}
}
\Bigl
[
\fr
{
(
{\bf t}_1
{\bf t}_1
)
}
{
(
{\bf t}_1
{\bf t}_2
)
}
\lf(
\sqrt
{h
(
{\bf t}_1
{\bf t}_2
)
\cos\al
+
u({\bf t}_1,{\bf t}_2)
\sin\al
}
\rg.
$$
\ses
\ses
\ses
\ses
$$
\lf. - \sqrt {\fr1h ( {\bf t}_1 {\bf t}_2 ) \cos\al + u({\bf
t}_1,{\bf t}_2) \sin\al } \rg) t_2^R + \sqrt {\fr1h ( {\bf t}_1
{\bf t}_2 ) \cos\al + u({\bf t}_1,{\bf t}_2) \sin\al } \,\, t_1^R
\Bigr]
$$
\ses
and
\ses\\
$$
f^R_p (g; {\bf t}_1, {\bf t}_2 ) t_2^p = \fr 1 { \sqrt { h \sqrt{
( {\bf t}_1 {\bf t}_1 )} \sqrt{ ( {\bf t}_2 {\bf t}_2 )} } } \sqrt
{h ( {\bf t}_1 {\bf t}_2 ) \cos\al + u({\bf t}_1,{\bf t}_2)
\sin\al } \,\, t_2^R,
$$
together with
\ses\\
$$
\sum_{R=1}^{N} f^R_p (g; {\bf t}_1, {\bf t}_2 ) t_1^R = \fr 1 {
\sqrt { h \sqrt{ ( {\bf t}_1 {\bf t}_1 )} \sqrt{ ( {\bf t}_2 {\bf
t}_2 )} } } \sqrt {h ( {\bf t}_1 {\bf t}_2 ) \cos\al + u({\bf
t}_1,{\bf t}_2) \sin\al } \,\, t_{1p}
$$
\ses\\
and
\ses\\
$$
\sum_{R=1}^{N}
f^R_p
(g;
{\bf t}_1,
{\bf t}_2
)
t_2^R
=
$$
\ses
\ses
$$
=
\fr
1
{
\sqrt
{
h
\sqrt{
(
{\bf t}_1
{\bf t}_1
)}
\sqrt{
(
{\bf t}_1
{\bf t}_2
)}
}
}
\Bigl
[
\fr
{
(
{\bf t}_2
{\bf t}_2
)
}
{
(
{\bf t}_1
{\bf t}_2
)
}
\lf(
\sqrt
{h
(
{\bf t}_1
{\bf t}_2
)
\cos\al
+
u({\bf t}_1,{\bf t}_2)
\sin\al
}
\rg.
$$
\ses
\ses
\ses
\ses
$$
\lf.
-
\sqrt
{\fr1h
(
{\bf t}_1
{\bf t}_2
)
\cos\al
+
u({\bf t}_1,{\bf t}_2)
\sin\al
}
\rg)
t_{1p}
+
\sqrt
{\fr1h
(
{\bf t}_1
{\bf t}_2
)
\cos\al
+
u({\bf t}_1,{\bf t}_2)
\sin\al
}
\,
\,t_{2p}
\Bigr].
$$

\ses
\ses
\ses
\ses

\setcounter{sctn}{3}
\setcounter{equation}{0}
3.3. {\it
Covariant version}.
It proves possible to convert the approach into the
{\it  co-version}
 by introducing
the {\it  co-vectors}
\be
 T_{1p}
(g;
{\bf t}_1,
{\bf t}_2
)
:=
n_{pq}
(g;
{\bf t}_1,
{\bf t}_2
)
{\bf t}_2^q,
\qquad
 T_{2q}
(g;
{\bf t}_1,
{\bf t}_2
)
:=
{\bf t}_1^p
n_{pq}(g;
{\bf t}_1,
{\bf t}_2
).
\ee
Applying (2.2),
we get
\ses\\
\be
{\bf T}_1
=
\fr
{
\sqrt
{
(
{\bf t}_2
{\bf t}_2
)
}
}
{
\sqrt
{
(
{\bf t}_1
{\bf t}_1
)
}
}
{\bf t}_1
\,
\cos
\al
+
\fr{
\sqrt{(
{\bf t}_2
{\bf t}_2
)}
}
{
h
\sqrt{(
{\bf t}_1
{\bf t}_1
)}
}
{\bf d}_1
\sin
\al
\ee
and
\be
{\bf T}_2
=
\fr
{
\sqrt
{
(
{\bf t}_1
{\bf t}_1
)
}
}
{
\sqrt
{
(
{\bf t}_2
{\bf t}_2
)
}
}
{\bf t}_2
\,
\cos
\al
+
\fr{
\sqrt{(
{\bf t}_1
{\bf t}_1
)}
}
{
h
\sqrt{(
{\bf t}_2
{\bf t}_2
)}
}
{\bf d}_2
\sin
\al.
\ee
The equality
\ses\\
\be
(
{\bf t}_1
{\bf T}_1
)
+
(
{\bf t}_2
{\bf T}_2
)
=
2
<{\bf t}_1,{\bf t}_2>
=
2
\sqrt{(
{\bf t}_1
{\bf t}_1
)}
\sqrt{(
{\bf t}_2
{\bf t}_2
)}
\cos\al
\ee
holds.
Also,
\be
\lim_{
{\bf t}_2\to
{\bf t}_1=
{\bf t}
}
\Bigl\{
{\bf T}_1
\Bigl\}=
\lim_{
{\bf t}_2\to
{\bf t}_1=
{\bf t}
}
\Bigl\{
{\bf T}_2
\Bigl\}=
{\bf t}.
\ee
The metric tensor (2.1)-(2.2) is obtainable from these vectors as follows:
\be
n_{pq}(g;
{\bf t}_1,
{\bf t}_2
)
=
\D
{T_{1p}}
{t_2^q}
=
\D
{T_{2q}}
{t_1^p}.
\ee

\ses

The respective products are found to be
\be
(
{\bf T}_1
{\bf T}_1
)
=
(
{\bf t}_2
{\bf t}_2
)
(\cos^2\al+\fr1{h^2}\sin^2\al),
\qquad
(
{\bf T}_2
{\bf T}_2
)
=
(
{\bf t}_1
{\bf t}_1
)
(\cos^2\al+\fr1{h^2}\sin^2\al),
\ee
and
\ses\\
\be
(
{\bf T}_1
{\bf T}_2
)
=
(\cos^2\al-\fr1{h^2}\sin^2\al)
(
{\bf t}_1
{\bf t}_2
)
+
2\fr1h
u({\bf t}_1,{\bf t}_2)
\cos
\al
\sin\al,
\ee
\ses
\ses\\
together with
\be
u({\bf T}_1,{\bf T}_2)
=
\fr2h
(
{\bf t}_1
{\bf t}_2
)
\sin\al\cos\al
-
(\cos^2\al-\fr1{h^2}\sin^2\al)
u({\bf t}_1,{\bf t}_2),
\ee
where
\be
u({\bf T}_1,{\bf T}_2)
=
\sqrt
{
(
{\bf T}_1
{\bf T}_1
)
(
{\bf T}_2
{\bf T}_2
)
-
(
{\bf T}_1
{\bf T}_2
)^2
}.
\ee

From  (3.7)-(3.8) it follows that
\ses\\
$$
(\cos^2\al+\fr1{h^2}\sin^2\al)^2
u({\bf t}_1,{\bf t}_2)
=
\fr2h
(
{\bf T}_1
{\bf T}_2
)
\sin\al\cos\al
-
(\cos^2\al-\fr1{h^2}\sin^2\al)
u({\bf T}_1,{\bf T}_2)
$$
and
\ses\\
$$
(\cos^2\al+\fr1{h^2}\sin^2\al)^2 ( {\bf t}_1 {\bf t}_2 ) =
(\cos^2\al-\fr1{h^2}\sin^2\al) ( {\bf T}_1 {\bf T}_2 ) + \fr2h
\sin\al\cos\al \,
u({\bf T}_1,{\bf T}_2),
$$
\ses
\ses
together with
\ses\\
$$
(\cos^2\al+\fr1{h^2}\sin^2\al) \Bigl[ -\fr1h ( {\bf t}_1 {\bf t}_2
) \sin\al + u({\bf t}_1,{\bf t}_2) \cos\al \Bigr] = \fr1h ( {\bf
T}_1 {\bf T}_2 ) \sin\al - u({\bf T}_1,{\bf T}_2) \cos\al.
$$

Using these formulas in calculating the co-representation
\be
\al=\hat\al( {\bf T}_1, {\bf T}_2 )
 \ee
 for the angle (1.21) yields
the following implicit equation:
\be \cos(h\al)= \fr {
(\cos^2\al-\fr1{h^2}\sin^2\al) ( {\bf T}_1 {\bf T}_2 ) + \fr2h
\sin\al\cos\al\, u({\bf T}_1,{\bf T}_2), } {
(\cos^2\al+\fr1{h^2}\sin^2\al) \sqrt{ ( {\bf T}_1 {\bf T}_1 ) }
\sqrt{ ( {\bf T}_2 {\bf T}_2 ) } }.
\ee

The respective co-version of the scalar product (1.36) reads
\be
<{\bf T}_1,{\bf T}_2>
=
\sqrt{
(
{\bf T}_1
{\bf T}_1
)
}
\sqrt{
(
{\bf T}_2
{\bf T}_2
)
}
\cos\al.
\ee

On this way the set (3.2)-(3.3) can be inverted, yielding
\ses\\
$$
{\bf t}_1
(g;
{\bf T}_1,
{\bf T}_2
)
=
\fr1{\xi}
\Bigl[
\fr
{
\sqrt{
(
{\bf t}_1
{\bf t}_1
)
}
}{
\sqrt{
(
{\bf t}_2
{\bf t}_2
)
}
}
\lf(
\cos\al-\fr1h
\fr
{
(
{\bf t}_1
{\bf t}_2
)
}
{
u({\bf t}_1,{\bf t}_2)
}
\sin\al
\rg)
{\bf T}_1
$$
\ses
\ses
\be
-
\fr1h
\fr
{
\sqrt{
(
{\bf t}_1
{\bf t}_1
)
}
\sqrt{
(
{\bf t}_2
{\bf t}_2
)
}
}
{
u({\bf t}_1,{\bf t}_2)
}
\sin\al
\,
{\bf T}_2
\Bigl]
\ee
and
\ses
\ses\\
$$
{\bf t}_2
(g;
{\bf T}_1,
{\bf T}_2
)
=
\fr1{\xi}
\Bigl[
\fr
{
\sqrt{
(
{\bf t}_2
{\bf t}_2
)
}
}{
\sqrt{
(
{\bf t}_1
{\bf t}_1
)
}
}
\lf(
\cos\al-\fr1h
\fr
{
(
{\bf t}_1
{\bf t}_2
)
}
{
u({\bf t}_1,{\bf t}_2)
}
\sin\al
\rg)
{\bf T}_2
$$
\ses
\ses
\be
-
\fr1h
\fr
{
\sqrt{
(
{\bf t}_1
{\bf t}_1
)
}
\sqrt{
(
{\bf t}_2
{\bf t}_2
)
}
}
{
u({\bf t}_1,{\bf t}_2)
}
\sin\al
\,\,
{\bf T}_1
\Bigl],
\ee
where
\ses\\
$$
{\xi}=
\lf(
\cos\al-\fr1h
\fr
{
(
{\bf t}_1
{\bf t}_2
)
}
{
u({\bf t}_1,{\bf t}_2)
}
\sin\al
\rg)^2
-
\fr1{h^2}
\lf
(
\fr
{
\sin\al
}
{
u({\bf t}_1,{\bf t}_2)
}
\rg)^2
(
{\bf t}_1
{\bf t}_1
)
(
{\bf t}_2
{\bf t}_2
),
$$
or
\ses\\
\be
{\xi}=
\cos^2\al-\fr1{h^2}\sin^2\al
-\fr2h
\fr
{
\sin\al
\cos\al
}
{
u({\bf t}_1,{\bf t}_2)
}
(
{\bf t}_1
{\bf t}_2
).
\ee
On taking into account (3.9), this function can be written merely as
\ses\\
\be
{\xi}=
-
\fr{
u({\bf T}_1,{\bf T}_2)
}
{
u({\bf t}_1,{\bf t}_2)
}.
\ee

Thus we find
\ses\\
\be
{\bf t}_1
=
\fr
{
\sqrt
{
(
{\bf T}_2
{\bf T}_2
)
}
}
{
\sqrt
{
(
{\bf T}_1
{\bf T}_1
)
}
}
{\bf T}_1
\,
\fr{\cos\al}
{\cos^2\al+\fr1{h^2}\sin^2\al}
+
\fr1h
\fr
{
\sqrt
{
(
{\bf T}_2
{\bf T}_2
)
}
}
{
\sqrt
{
(
{\bf T}_1
{\bf T}_1
)
}
}
{\bf D}_1
\fr{\sin\al}
{\cos^2\al+\fr1{h^2}\sin^2\al}
\ee
and
\be
{\bf t}_2
=
\fr
{
\sqrt
{
(
{\bf T}_1
{\bf T}_1
)
}
}
{
\sqrt
{
(
{\bf T}_2
{\bf T}_2
)
}
}
{\bf T}_2
\,
\fr{\cos\al}
{\cos^2\al+\fr1{h^2}\sin^2\al}
+
\fr1h
\fr
{
\sqrt
{
(
{\bf T}_1
{\bf T}_1
)
}
}
{
\sqrt
{
(
{\bf T}_2
{\bf T}_2
)
}
}
{\bf D}_2
\fr{\sin\al}
{\cos^2\al+\fr1{h^2}\sin^2\al},
\ee
where
\ses\\
\be
{\bf D}_1
=
\fr{
(
{\bf T}_1
{\bf T}_1
)
{\bf T}_2
-
(
{\bf T}_1
{\bf T}_2
)
{\bf T}_1
}
{
u({\bf T}_1,{\bf T}_2)
},
\qquad
{\bf D}_2
=
\fr{
(
{\bf T}_2
{\bf T}_2
)
{\bf T}_1
-
(
{\bf T}_1
{\bf T}_2
)
{\bf T}_2
}
{
u({\bf T}_1,{\bf T}_2)
}.
\ee

\ses
\ses
\ses

The identities
$$
(
{\bf T}_1
{\bf D}_1
)=0,
\qquad
(
{\bf T}_2
{\bf D}_2
)=0,
$$
\ses
$$
(
{\bf D}_1
{\bf D}_2
)=-
(
{\bf T}_1
{\bf T}_2
),
\qquad
(
{\bf D}_1
{\bf D}_1
)=
(
{\bf T}_1
{\bf T}_1
),
\qquad
(
{\bf D}_2
{\bf D}_2
)=
(
{\bf T}_2
{\bf T}_2
),
$$
and
\ses\\
$$
(
{\bf D}_1
{\bf T}_2
)
=
(
{\bf T}_1
{\bf D}_2
)
=
u({\bf T}_1,{\bf T}_2)
$$
hold.

By the help of (3.20)-(3.21), and in close similarity to (3.6),
the co-version
$$
N^{pq}(g;
{\bf T}_1,
{\bf T}_2
)
:=
\D
{t_1^p}
{T_{2q}}
=
\D
{t_2^q}
{T_{1p}}
$$
for the two-vector metric tensor (2.2) can be arrived at.

\ses
\ses

\setcounter{sctn}{4}
\setcounter{equation}{0}
3.4. {\it
$\cE_g^{PD}$-parallelogram law}.
Let
$
{\bf t}_1,
 {\bf t}_2$, and
$
{\bf t}_3
$
be three vectors issued from the same origin $``O"$,
subject to the conditions that the angle between
$
{\bf t}_1
$
and
$
{\bf t}_2
$
is acute and the vector
$
{\bf t}_3
$
is positioned between the vectors
$
{\bf t}_1
$
and
$
{\bf t}_2
$.
Let us denote the end points of the vectors
$
{\bf t}_1,
 {\bf t}_2$, and
$
{\bf t}_3
$
as
$X_1, X_2$, and $X_3$,
respectively.
On joining the points
$X_1$ and $X_3$, and also
$X_2$ and $X_3$,
by  means of
$\cE_g^{PD}$-geodesics, we get a tetragonal figure, to be denoted as
$\cP_4$.

Using Eqs. (1.16), (1.17), and (1.35), we  set forth
the following couple  equations:
\be
({\bf t}_2
{\bf t}_2
)
=
({\bf t}_1
{\bf t}_1
)
+
({\bf t}_3
{\bf t}_3
)
-2
\sqrt
{
({\bf t}_1
{\bf t}_1
)
}
\,
\sqrt
{
({\bf t}_3
{\bf t}_3
)
}
\cos
\Bigl[
\fr1h\arccos
\fr{
(
{\bf t}_1
{\bf t}_3
)}
{
\sqrt{(
{\bf t}_1
{\bf t}_1
)}
\,
\sqrt{(
{\bf t}_3
{\bf t}_3
)}
}
\Bigr]
\ee
and
\ses\\
\be
({\bf t}_1
{\bf t}_1
)
=
({\bf t}_3
{\bf t}_3
)
+
({\bf t}_2
{\bf t}_2
)
-2
\sqrt
{
({\bf t}_3
{\bf t}_3
)
}
\,
\sqrt
{
({\bf t}_2
{\bf t}_2
)
}
\cos
\Bigl[
\fr1h\arccos
\fr{
(
{\bf t}_2
{\bf t}_3
)}
{
\sqrt{(
{\bf t}_2
{\bf t}_2
)}
\,
\sqrt{(
{\bf t}_3
{\bf t}_3
)}
}
\Bigr],
\ee
\ses\\
which can also be rewritten in the convenient form
\be
\sqrt{(
{\bf t}_3
{\bf t}_3
)}
-
\fr{
({\bf t}_2
{\bf t}_2
)
-
({\bf t}_1
{\bf t}_1
)
}
{
\sqrt{(
{\bf t}_3
{\bf t}_3
)}
}
=2
\sqrt
{
({\bf t}_1
{\bf t}_1
)
}
\cos
\Bigl[
\fr1h\arccos
\fr{
(
{\bf t}_1
{\bf t}_3
)}
{
\sqrt{(
{\bf t}_1
{\bf t}_1
)}
\,
\sqrt{(
{\bf t}_3
{\bf t}_3
)}
}
\Bigr]
\ee
and
\ses\\
\be
\sqrt{(
{\bf t}_3
{\bf t}_3
)}
-
\fr{
({\bf t}_1
{\bf t}_1
)
-
({\bf t}_2
{\bf t}_2
)
}
{
\sqrt{(
{\bf t}_3
{\bf t}_3
)}
}
=2
\sqrt
{
({\bf t}_2
{\bf t}_2
)
}
\cos
\Bigl[
\fr1h\arccos
\fr{
(
{\bf t}_2
{\bf t}_3
)}
{
\sqrt{(
{\bf t}_2
{\bf t}_2
)}
\,
\sqrt{(
{\bf t}_3
{\bf t}_3
)}
}
\Bigr].
\ee

In (4.1), the left-hand part is the squared length of the straight side $OX_2$
and the
right-hand side is the squared length of the geodesic side
$X_1X_3$. According to (4.2),
the lengths of $OX_1$
and $X_2X_3$ are equal.
Under these conditions, the figure $\cP_4$ does attribute the general
property of the Euclidean parallelogram that the lengths of opposite sides
are equal. In this vein, we introduce the following

\ses
\ses

{\large  Definition}. Subject to the equations (4.1) and (4.2),
the tetragonal figure
$\cP_4$
is called
the
{
\it
$\cE_g^{PD}$-parallelogram,}
and the vector
$
{\bf t}_3
$
is called
the
$\cE_g^{PD}$-sum vector:
\be
{\bf t}_3
=
{\bf t}_1
\oplus
{\bf t}_2.
\ee

{\large Note}. The qualitative distinction here
from Euclidean patterns
is that the  sides
$X_1X_3$ and $X_3X_2$ of the $\cP_4$ are curved lines in general,
namely  geodesic arcs,
which generally cease to be straight under the $\cE_g^{PD}$-extension.

\ses\ses

\begin{figure}[t]
\centering
\includegraphics[width=8cm]{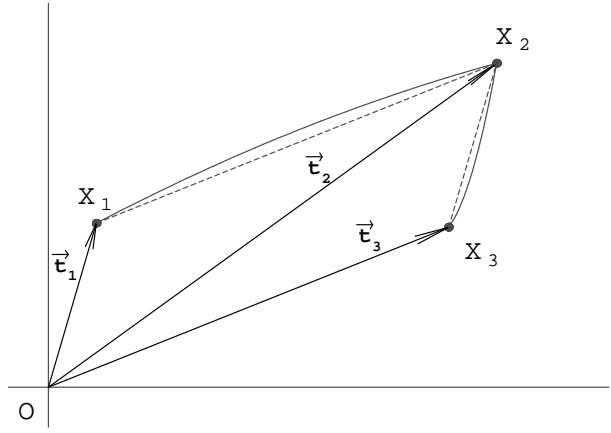}
\caption{ [The Finsleroid-parallelogram law is applied]}
\end{figure}

\ses\ses
\ses\ses

Finding the sum vector (4.5)
implies solving the set of the equations (4.3) and (4.4).
We shall proceed  approximately, namely taking
\be
\fr1h=1+k
\ee
and
\be
{\bf t}_3
=
{\bf t}_1
+{\bf t}_2
+k
{\bf c}
({\bf t}_1,{\bf t}_2),
\quad
k\ll 1.
\ee
Under these conditions, on inserting (4.6) and (4.7) in (4.3), we find
\ses\\
$$
\sqrt
{
(
{\bf t}_2
+{\bf t}_1
)^2
}
+
k\fr{
(
{\bf t}_2
+{\bf t}_1
){\bf c}
}
{
\sqrt
{
(
{\bf t}_2
+{\bf t}_1
)^2
}
}
-
\fr{
({\bf t}_2
{\bf t}_2
)
-
({\bf t}_1
{\bf t}_1
)
}
{
\sqrt
{
(
{\bf t}_2
+{\bf t}_1
)^2
}
}
\lf(1-k\fr{
(
{\bf t}_2
+{\bf t}_1
){\bf c}
}
{
{
(
{\bf t}_2
+{\bf t}_1
)^2
}
}
\rg)
$$
\ses
\ses
$$
=2
\sqrt
{
({\bf t}_1
{\bf t}_1
)
}
\cos
\Bigl[
(1+k)\arccos
\fr{
(
{\bf t}_1
{\bf t}_3
)}
{
\sqrt{(
{\bf t}_1
{\bf t}_1
)}
\,
\sqrt{(
{\bf t}_3
{\bf t}_3
)}
}
\Bigr]
$$
\ses
\ses
\ses
\ses
$$
=2
\fr{
(
{\bf t}_1
{\bf t}_3
)}
{
\sqrt{(
{\bf t}_3
{\bf t}_3
)}
}
-
2
k
\sqrt
{
({\bf t}_1
{\bf t}_1
)
}
\sqrt
{
1-
\lf(
\fr{
{\bf t}_1
(
{\bf t}_2
+
{\bf t}_1
)
}
{
\sqrt{(
{\bf t}_1
{\bf t}_1
)}
\,
\sqrt
{
(
{\bf t}_2
+{\bf t}_1
)^2
}
}
\rg)^2
}
\arccos
\fr{
{\bf t}_1
(
{\bf t}_2
+
{\bf t}_1
)
}
{
\sqrt{(
{\bf t}_1
{\bf t}_1
)}
\,
\sqrt
{
(
{\bf t}_2
+{\bf t}_1
)^2
}
}
$$
\ses
\ses
\ses
\ses
\ses
$$
=
2k
\fr{
(
{\bf t}_1
{\bf c}
)}
{
\sqrt
{
(
{\bf t}_2
+{\bf t}_1
)^2
}
}
+
2
\fr{
{\bf t}_1
(
{\bf t}_2
+
{\bf t}_1
)}
{
\sqrt
{
(
{\bf t}_2
+{\bf t}_1
)^2
}
}
\lf(1-k
\fr{
(
{\bf t}_2
+{\bf t}_1
){\bf c}
}
{
{
(
{\bf t}_2
+{\bf t}_1
)^2
}
}
\rg)
$$
\ses
\ses
$$
-
2k
\sqrt
{
({\bf t}_1
{\bf t}_1
)
}
\sqrt
{
1-
\lf(
\fr{
{\bf t}_1
(
{\bf t}_2
+
{\bf t}_1
)
}
{
\sqrt{(
{\bf t}_1
{\bf t}_1
)}
\,
\sqrt
{
(
{\bf t}_2
+{\bf t}_1
)^2
}
}
\rg)^2
}
\arccos
\fr{
{\bf t}_1
(
{\bf t}_2
+
{\bf t}_1
)
}
{
\sqrt{(
{\bf t}_1
{\bf t}_1
)}
\,
\sqrt
{
(
{\bf t}_2
+{\bf t}_1
)^2
}
},
$$
\ses

\clearpage

which entails
$$
\fr{
(
{\bf t}_2
+{\bf t}_1
){\bf c}
}
{
\sqrt
{
(
{\bf t}_2
+{\bf t}_1
)^2
}
}
+
\fr{
({\bf t}_2
{\bf t}_2
)
-
({\bf t}_1
{\bf t}_1
)
}
{
{
(
{\bf t}_2
+{\bf t}_1
)^2
}
}
\fr{
(
{\bf t}_2
+{\bf t}_1
){\bf c}
}
{
\sqrt
{
(
{\bf t}_2
+{\bf t}_1
)^2
}
}
$$
\ses
\ses
\ses
\ses
\ses
$$
=
2
\fr{
(
{\bf t}_1
{\bf c}
)}
{
\sqrt
{
(
{\bf t}_2
+{\bf t}_1
)^2
}
}
-
2
\fr{
{\bf t}_1
(
{\bf t}_2
+
{\bf t}_1
)}
{
{
(
{\bf t}_2
+{\bf t}_1
)^2
}
}
\fr{
(
{\bf t}_2
+{\bf t}_1
){\bf c}
}
{
\sqrt
{
(
{\bf t}_2
+{\bf t}_1
)^2
}
}
$$
\ses
\ses
$$
-
2
\sqrt
{
({\bf t}_1
{\bf t}_1
)
}
\sqrt
{
1-
\lf(
\fr{
{\bf t}_1
(
{\bf t}_2
+
{\bf t}_1
)
}
{
\sqrt{(
{\bf t}_1
{\bf t}_1
)}
\,
\sqrt
{
(
{\bf t}_2
+{\bf t}_1
)^2
}
}
\rg)^2
}
\arccos
\fr{
{\bf t}_1
(
{\bf t}_2
+
{\bf t}_1
)
}
{
\sqrt{(
{\bf t}_1
{\bf t}_1
)}
\,
\sqrt
{
(
{\bf t}_2
+{\bf t}_1
)^2
}
}.
$$

Therefore we obtain
$$
{\bf t}_2
{\bf c}
=
-
u
(
{\bf t}_1,
{\bf t}_2
)
\arccos
\fr{
{\bf t}_1
(
{\bf t}_2
+
{\bf t}_1
)
}
{
\sqrt{(
{\bf t}_1
{\bf t}_1
)}
\,
\sqrt
{
(
{\bf t}_2
+{\bf t}_1
)^2
}
},
$$
\ses where $ u ( {\bf t}_1, {\bf t}_2 ) $ is the function (1.20).
Similarly, from (4.4) it follows that
\ses\\
$$
{\bf t}_1
{\bf c}
=
-
u
(
{\bf t}_1,
{\bf t}_2
)
\arccos
\fr{
{\bf t}_2
(
{\bf t}_2
+
{\bf t}_1
)
}
{
\sqrt{(
{\bf t}_2
{\bf t}_2
)}
\,
\sqrt
{
(
{\bf t}_2
+
{\bf t}_1
)^2
}
}.
$$

If we use now the symmetrized expansion
$$
{\bf c} =
m
(
{\bf t}_1,
{\bf t}_2
)
{\bf t}_1
+
n
(
{\bf t}_1 ,
{\bf t}_2
)
{\bf t}_2,
$$
then we find
\ses\\
$$
m
(
{\bf t}_1,
{\bf t}_2
)
=
$$
\ses
\ses
$$
\fr1
{
u
(
{\bf t}_1,
{\bf t}_2
)
}
\lf(
(
{\bf t}_1
{\bf t}_2
)
\arccos
\fr{
{\bf t}_1
(
{\bf t}_2
+
{\bf t}_1
)
}
{
\sqrt{(
{\bf t}_1
{\bf t}_1
)}
\,
\sqrt
{
(
{\bf t}_2
+{\bf t}_1
)^2
}
}
-
(
{\bf t}_2
{\bf t}_2
)
\arccos
\fr{
{\bf t}_2
(
{\bf t}_2
+
{\bf t}_1
)
}
{
\sqrt{(
{\bf t}_2
{\bf t}_2
)}
\,
\sqrt
{
(
{\bf t}_2
+{\bf t}_1
)^2
}
}
\rg)
$$
\ses
and
\ses\\
$$
n
(
{\bf t}_1,
{\bf t}_2
)
=
$$
\ses
\ses
$$
\fr1
{
u
(
{\bf t}_1,
{\bf t}_2
)
}
\lf(
(
{\bf t}_1
{\bf t}_2
)
\arccos
\fr{
{\bf t}_2
(
{\bf t}_2
+
{\bf t}_1
)
}
{
\sqrt{(
{\bf t}_2
{\bf t}_2
)}
\,
\sqrt
{
(
{\bf t}_2
+{\bf t}_1
)^2
}
}
-
(
{\bf t}_1
{\bf t}_1
)
\arccos
\fr{
{\bf t}_1
(
{\bf t}_2
+
{\bf t}_1
)
}
{
\sqrt{(
{\bf t}_1
{\bf t}_1
)}
\,
\sqrt
{
(
{\bf t}_2
+{\bf t}_1
)^2
}
}
\rg).
$$

Since
$$
m
(
{\bf t}_1,
{\bf t}_2
)
=
n
(
{\bf t}_2,
{\bf t}_1
),
$$
 we  just deduce the approximate solution
\ses
$$
{\bf t}_1
\oplus
{\bf t}_2
\approx
{\bf t}_1
+{\bf t}_2
+(\fr1h-1)
\Bigl
(
m
({\bf t}_1,{\bf t}_2)
{\bf t}_1
+
m
({\bf t}_2,{\bf t}_1)
{\bf t}_2
\Bigr
),
\qquad
\fr1h-1\ll 1.
$$

Alternatively, the solution
$
{\bf t}_2
=
{\bf t}_2
({\bf t}_1,{\bf t}_3)
$
to the set of equations (4.1)-(4.2) can naturally be treated as
the
{
\it
$\cE_g^{PD}$-difference} of vectors
$
{\bf t}_3
$
and
$
{\bf t}_1
$:
$$
{\bf t}_2
=
{\bf t}_3
\ominus
{\bf t}_1.
$$
Again, restricting ourselves to the approximation,
 from (4.1)-(4.2) we obtain
$$
{\bf t}_3
\ominus
{\bf t}_1
\approx
{\bf t}_3
-
{\bf t}_1
+
(
\fr1h-1
)
{\bf s}
(
{\bf t}_1,
{\bf t}_3
), \qquad
\fr1h-1\ll 1,
$$
with
$$
{\bf s}
(
{\bf t}_1,
{\bf t}_3
)
=
\fr1
{
u
(
{\bf t}_1,
{\bf t}_3
)
}
\Bigl\{
\Bigl[
(
{\bf t}_1
{\bf t}_1
)
\arccos
\fr{
(
{\bf t}_1
{\bf t}_3
)
}
{
\sqrt{(
{\bf t}_1
{\bf t}_1
)}
\,
\sqrt
{
(
{\bf t}_3
{\bf t}_3
)
}
}
$$
\ses
\ses
$$
-
(
{\bf t}_3
-
{\bf t}_1,
{\bf t}_1
)
\arccos
\fr{
(
{\bf t}_3-
{\bf t}_1,
{\bf t}_3
)
}
{
\sqrt{(
{\bf t}_3 -
{\bf t}_1,
{\bf t}_3 -
{\bf t}_1
)}
\,
\sqrt
{
(
{\bf t}_3
{\bf t}_3
)
}
}
\Bigr]
(
{\bf t}_3-
{\bf t}_1
)
$$
\ses
\ses
$$
+
\Bigl[
(
{\bf t}_3-
{\bf t}_1,
{\bf t}_3 -
{\bf t}_1
)
\arccos
\fr{
(
{\bf t}_3-
{\bf t}_1,
{\bf t}_3
)
}
{
\sqrt{(
{\bf t}_3 -
{\bf t}_1,
{\bf t}_3 -
{\bf t}_1
)}
\,
\sqrt
{
(
{\bf t}_3
{\bf t}_3
)
}
}
$$
\ses
\ses
$$
-
(
{\bf t}_3
-
{\bf t}_1,
{\bf t}_1
)
\arccos
\fr{
(
{\bf t}_1
{\bf t}_3
)
}
{
\sqrt{(
{\bf t}_1
{\bf t}_1
)}
\,
\sqrt
{
(
{\bf t}_3
{\bf t}_3
)
}
}
\Bigr]
{\bf t}_1
\Bigl\}.
$$
\ses
Here it is useful to note that
$$
(
{\bf t}_3
-
{\bf t}_1,
{\bf s}
)
= u
(
{\bf t}_1,
{\bf t}_3
)
\arccos
\fr{
(
{\bf t}_1
{\bf t}_3
)
}
{
\sqrt{(
{\bf t}_1
{\bf t}_1
)}
\,
\sqrt
{
(
{\bf t}_3
{\bf t}_3
)
}
},
$$
\ses
\ses
$$
(
{\bf t}_1,
{\bf s}
)
=
u
(
{\bf t}_1,
{\bf t}_3
)
\arccos
\fr{
(
{\bf t}_3-
{\bf t}_1,
{\bf t}_3
)
}
{
\sqrt{(
{\bf t}_3-
{\bf t}_1,
{\bf t}_3-
{\bf t}_1
)}
\,
\sqrt
{
(
{\bf t}_3
{\bf t}_3
)
}
},
$$
\ses\\
and
$$
u
(
{\bf t}_3-
{\bf t}_1,
{\bf t}_3
)
=
u
(
{\bf t}_1,
{\bf t}_3
).
$$
The problem of finding
the relevant sum vector
$
{\bf t}_1
\oplus
{\bf t}_2
$
and
 difference vector
$
{\bf t}_3
\ominus
{\bf t}_2
$
in  general exact forms
is open and seems to be difficult.

\bigskip
\bigskip

\def\bibit[#1]#2\par{\rm\noindent\parskip1pt
                     \parbox[t]{.05\textwidth}{\mbox{}\hfill[#1]}\hfill
                     \parbox[t]{.925\textwidth}{\baselineskip11pt#2}\par}

\bc
 References
\ec
\bigskip

\bibit[1] E. Cartan: \it Les espaces de Finsler, Actualites \rm 79,
Hermann,
Paris 1934.

\bibit[2] H. Busemann: \it Canad. J. Math. \bf1 \rm(1949), 279.

\bibit[3] H. Rund: \it The Differential Geometry of Finsler spaces, \rm
Springer-Verlag, Berlin 1959.

\bibit[4] G. S. Asanov: \it Finsler Geometry, Relativity and Gauge Theories, \rm
D.~Reidel Publ. Comp., Dordrecht 1985.

\bibit[5] R. S. Ingarden: \it Tensor \bf30 \rm(1976), 201.

\bibit[6] R. S. Ingarden and L. Tamassy: \it  \it Rep. Math. Phys. \bf 32 \rm(1993), 11.

\bibit[7] I. Kozma and L. Tamassy: \it  \it Rep. Math. Phys. \bf 42 \rm(2003), 77.

\bibit[8] D.~Bao, S. S. Chern, and Z. Shen: \it An
Introduction to Riemann-Finsler Geometry,
\quad\rm Springer, N.Y., Berlin, 2000.

\bibit[9] A. C. Thompson: \it Minkowski Geometry \rm (Encycl. of Math. and its Appl.
{\bf 63}), Cambridge Univ. Press, Cambridge. 1996.

\bibit[10] G. S. Asanov: \it Aeq. Math. \bf49 \rm(1995), 234.

\bibit[11] G.S. Asanov: \it Rep. Math. Phys. \bf 45 \rm(2000), 155;
\bf 47 \rm(2001), 323.

\bibit[12] G. S. Asanov: arXiv:math-ph/0310019, 2003.

\end{document}